\documentclass[11pt]{amsart}
\usepackage[T1]{fontenc}
\usepackage{lmodern,amsmath,amssymb,mathtools,microtype}
\usepackage{booktabs}
\usepackage{cite}
\usepackage{enumitem}
\usepackage[letterpaper,margin=1.05in]{geometry}
\usepackage[colorlinks=true,linkcolor=blue,citecolor=blue,urlcolor=blue]{hyperref}
\numberwithin{equation}{section}
\newtheorem{theorem}{Theorem}[section]
\newtheorem{lemma}[theorem]{Lemma}
\newtheorem{corollary}[theorem]{Corollary}
\newtheorem{proposition}[theorem]{Proposition}
\newcommand{\R}{\mathbb R}
\newcommand{\sphere}{\mathbb S^{N-1}}
\newcommand{\En}{\mathcal E_{s,p}}
\newcommand{\al}{\alpha_{s,N}}
\newcommand{\gm}{\gamma_{s,N}}

\newtheorem*{thmPR}{Theorem B}
\newtheorem*{thmM}{Theorem A}
\title[The sharp fractional Trudinger--Moser exponent]{The sharp Trudinger--Moser exponent in Sobolev--Slobodeckij spaces}

\author[Phi Long Le]{Phi Long Le\textsuperscript{a}}
\author[Zhenzhen Wei]{Zhenzhen Wei\textsuperscript{b,*}}

\makeatletter
\renewcommand{\@setauthors}{%
  \par
  \vspace{6pt}%
  \begingroup
    \centering
    \normalfont\normalsize
    \def\and{, }%
    \authors\par
  \endgroup
}
\makeatother

\thanks{\textsuperscript{a}\,School of Economic Mathematics and Statistics,
  CELG, University of Economics Ho Chi Minh City, Ho Chi Minh City, Vietnam.}

\thanks{\textsuperscript{b}\,Department of Mathematics and Statistics,
  Memorial University of Newfoundland, St. John's, Newfoundland and Labrador
  A1C 5S7, Canada.}

\thanks{\textit{E-mail addresses.}
  \texttt{phill@ueh.edu.vn} (Phi Long Le),
  \texttt{zhenzhenw@mun.ca} (Zhenzhen Wei).}
\thanks{\textsuperscript{*}\,Corresponding author.}
\date{}
\allowdisplaybreaks
\begin{document}
\begin{abstract}
We establish the sharp Trudinger--Moser inequality in the fractional Sobolev--Slobodeckij spaces $W^{s,p}$, $0 < s < 1$, at the critical index $sp=N$, including the endpoint. More precisely, there exists a constant $C_{N,s}$ such that for every measurable set $\Omega\subset\mathbb{R}^N$ of finite positive measure,
$$
\sup\Bigl\{\frac{1}{|\Omega|}\int_\Omega e^{\alpha_{s,N}|u|^{\frac{p}{p-1}}}dx: u\in W^{s,p}(\mathbb{R}^N), u=0\text{ a.e. on }\mathbb{R}^N\setminus\Omega, [u]_{s,p}\le 1\Bigr\}\le C_{N,s},
$$
where $[u]_{s,p}$ is the Gagliardo seminorm and $\alpha_{s,N}$ is the upper threshold found by Parini and Ruf for $N\ge 2$ and by Iula for $N=1$. Moreover, the constant $\alpha_{s,N}$ is sharp. This answers two of their questions.
We use a nonlocal version of Moser's one-dimensional lemma. The $p$-th moment of a translation kernel on the logarithmic half-line determines the sharp coefficient, and finiteness of the next moment characterizes the endpoint. We obtain sharp singular, whole-space, Adachi--Tanaka type, concentration-compactness and affine inequalities. When the Euclidean distance in the seminorm is replaced by the gauge of an origin-symmetric star body $K$, we determine the sharp exponent for each fixed $K$ and show that it exceeds the volume-normalized exponent unless $K$ is an ellipsoid. The same exponents are sharp for regional seminorms on bounded Lipschitz domains in dimensions $N\ge 2$.
We prove attainment of the critical supremum on balls for $N\ge 2$ and $s$ close to one. In the Hilbert case $N=1$, $s=1/2$, we compute the concentration level of the singular problem on an interval using a sharp weighted logarithmic Hardy--Littlewood--Sobolev inequality. Consequently, the critical supremum is attained for every singular weight $|x|^{\lambda-1}$, $0 < \lambda\le 1$.
\end{abstract}
\maketitle

\section{Introduction}
The main purpose of this paper is to establish the sharp Trudinger--Moser inequality in the fractional Sobolev--Slobodeckij spaces $W^{s,p}(\R^N)$ at the critical index $sp=N$, including the endpoint, and to study its singular, whole-space, affine, anisotropic and regional versions together with the existence of extremal functions.

We first recall the classical Trudinger--Moser inequality. Let $\Omega\subset\R^N$, $N\ge2$, be a bounded domain. The Sobolev embedding theorem tells us that the embedding $W_0^{1,q}(\Omega)\subset L^r(\Omega)$ is continuous for $1\le r\le Nq/(N-q)$ when $q<N$. However, $W_0^{1,N}(\Omega)\not\subset L^\infty(\Omega)$. In the borderline case, Yudovich \cite{Yudovich}, Pohozaev \cite{Pohozaev} and Trudinger \cite{Trudinger} proved that $W_0^{1,N}(\Omega)$ is embedded in the Orlicz space associated with the Young function $\exp(\beta|t|^{N/(N-1)})-1$ for some $\beta>0$. In 1971, Moser \cite{Moser} sharpened this result by finding the optimal $\beta$.
\begin{thmM}[Moser \cite{Moser}]
Let $\Omega\subset\R^N$, $N\ge2$, be a domain of finite measure. Then there exists a constant $c_0$ such that
\[
 \frac1{|\Omega|}\int_\Omega\exp\bigl(\beta_N|u|^{\frac N{N-1}}\bigr)dx\le c_0
\]
for all $u\in W_0^{1,N}(\Omega)$ with $\int_\Omega|\nabla u|^Ndx\le1$, where $\beta_N=N|\sphere|^{\frac1{N-1}}$. Moreover, the constant $\beta_N$ is sharp in the sense that if we replace $\beta_N$ by any $\beta>\beta_N$, then the above inequality can no longer hold with some $c_0$ independent of $u$.
\end{thmM}
Moser's proof reduces the problem by Schwarz symmetrization to radial functions, and then to a one-dimensional inequality for functions on a half-line; see also Jodeit \cite{Jodeit}. In 1988, Adams \cite{Adams} extended Theorem A to higher order Sobolev spaces by representing $u$ as a Riesz potential of its derivatives. Since then, the Trudinger--Moser and Adams inequalities have been extended in many directions, for instance to unbounded domains and to the whole space \cite{AT,Ruf,LiRuf,LamLu12a,LLZ,VHN19}, to singular weights \cite{AS,LamLu12c}, to the exact growth condition \cite{LuTangZhu15}, to the Heisenberg group and to groups of Heisenberg type \cite{CohnLu01,CohnLu02,LamLu12b,LiLuZhu18}, to hyperbolic spaces \cite{LuTang13,LuTang16,LiLuYang20,NgoNguyen20}, to homogeneous weights \cite{DuyNghiaLe19}, and to anisotropic, Finsler and affine energies \cite{LamLuTang17,DuyLamLe21,DuyLe22}. In several of these settings symmetrization is not available, and the proofs use instead representation formulas and fundamental solutions, as in the work of Cohn and Lu \cite{CohnLu01,CohnLu02}, or estimates on level sets, as in the work of Lam and Lu \cite{LamLu}. The method of the present paper is in this spirit. The existence of extremal functions on balls was proved by Carleson and Chang \cite{CC}, and their one-dimensional analysis was extended to real exponents by Hudson and Leckband \cite{HL}. For the existence and nonexistence of extremal functions in other settings we refer the interested reader to \cite{LuYang09,LLZ19a,VHN21a,ChenLuZhu22}, among others.

For fractional orders the situation is different. When the energy is the $L^p$ norm of a fractional Laplacian, a version of Adams' method applies, and sharp inequalities were obtained by Martinazzi \cite{Martinazzi}. We also mention the work of Takahashi \cite{Takahashi} on the whole line. In the Sobolev--Slobodeckij setting, Zhang \cite{Zhang} extended the inequality of \cite{PR} to the whole space. For every exponent below the optimal exponent of \cite{PR}, Chen, Lu and Zhang \cite{ChenLuZhang21} proved the existence of extremals for the subcritical fractional Trudinger--Moser inequalities on $\R^N$, by the level-set method of Lam and Lu and fractional rearrangement inequalities, and they also gave a concentration-compactness principle. Lu and Temperini \cite{LuTemperini25} established Trudinger--Moser inequalities in fractional Sobolev--Slobodeckij spaces on bounded domains of the Heisenberg group and on the whole group, for exponents below an optimal exponent whose value is not known, together with the corresponding concentration-compactness principles; for the global inequality they developed a technique which avoids symmetrization. In dimension one, Iula, Maalaoui and Martinazzi \cite{IMM} proved the sharp inequality with the constant $\pi$ for functions supported in a bounded interval $I$ under the constraint $\|(-\Delta)^{1/4}u\|_{L^2(I)}\le1$ \cite[Theorem~1.1]{IMM}, and on the whole line under the full norm $\|u\|_{L^2(\R)}^2+\|(-\Delta)^{1/4}u\|_{L^2(\R)}^2\le1$ \cite[Theorem~1.5]{IMM}; the homogeneous constraint alone does not suffice on $\R$, by dilation. Mancini and Martinazzi \cite{MM} proved the existence of extremals on intervals. However, the natural energy on $W^{s,p}(\R^N)$ is the Gagliardo seminorm, and the norm of $W^{s,p}(\R^N)$ is $(\|u\|_p^p+[u]_{s,p}^p)^{1/p}$; the seminorm alone does not control $\|u\|_p$ on $\R^N$, since it is invariant under dilations when $sp=N$. Throughout this paper we use
\begin{equation}\label{eq:seminorm}
 [u]_{s,p}^p=\int_{\R^N}\int_{\R^N}
 \frac{|u(x)-u(y)|^p}{|x-y|^{2N}}\,dx\,dy,\qquad p=\frac Ns .
\end{equation}
We refer to Di Nezza, Palatucci and Valdinoci \cite{DNPV} for the basic properties of the fractional Sobolev spaces. For $p\ne2$ this is not the $L^p$ norm of a fractional derivative, and the representation formula of Adams does not apply. For a measurable set \(\Omega\), let \(X(\Omega)\) denote the functions in \(W^{s,p}(\R^N)\) that vanish almost everywhere on its complement. The space \(X(\Omega)\) contains the closure of \(C_c^\infty(\Omega)\) in the zero-extension \(W^{s,p}\) norm, since convergence in that norm implies convergence in \(L^p(\R^N)\). The interactions of \(\Omega\) with its complement are part of \eqref{eq:seminorm}. Parini and Ruf \cite{PR} computed the Gagliardo energy of radial functions, and they proved the following result.
\begin{thmPR}[Parini and Ruf {\cite[Theorem 1.1]{PR}}]
Let $\Omega\subset\R^N$, $N\ge2$, be a bounded domain with Lipschitz boundary, let $0<s<1$, $sp=N$, and let $W_0^{s,p}(\Omega)$ be the completion of $C_c^\infty(\Omega)$ with respect to the norm $(\|u\|_{L^p(\Omega)}^p+[u]_{s,p}^p)^{1/p}$. Then there exists $\alpha^*=\alpha^*(s,\Omega)>0$ such that
\[
 \sup\Bigl\{\int_\Omega\exp\bigl(\alpha|u|^{\frac N{N-s}}\bigr)dx:\ u\in W_0^{s,p}(\Omega),\ [u]_{s,p}\le1\Bigr\}<\infty
\]
for $\alpha\in[0,\alpha^*)$. Moreover, the supremum is infinite for $\alpha\in(\alpha^*_{s,N},\infty)$, where $\alpha^*_{s,N}$ is an explicit constant.
\end{thmPR}
Following \cite[Section~6]{PR}, we call the supremum of all $\alpha$ for which the supremum in Theorem B is finite the optimal exponent, and we denote it by $\alpha^*_{\rm opt}(s,\Omega)$. By Theorem B, $0<\alpha^*_{\rm opt}(s,\Omega)\le\alpha^*_{s,N}$.
In dimension one Iula \cite[Theorems~1.1--1.3]{Iula} identified the corresponding upper threshold, proved finiteness for a positive range of coefficients, and established divergence at the threshold after multiplication by a factor tending to infinity. The authors of \cite{PR} asked whether $\alpha^*_{s,N}$ is the optimal exponent, whether the inequality holds at the optimal exponent, and whether a maximizer exists \cite[Section~6, p.~299]{PR}. It then becomes a very interesting and natural question to determine the sharp exponent, and to decide whether the inequality holds at the endpoint. In this paper we will show that the answer to the first two questions is positive in all dimensions, and we will give a partial answer to the third one.

To state our results, let us agree to some conventions. All functions may be taken real-valued. Write \(S_N=|\sphere|=N\omega_N\), where \(\omega_N=|B_1|\), with \(S_1=2\), and \(p'=p/(p-1)\). Define
\begin{equation}\label{eq:constant}
 \gm=\frac{2S_N^2\Gamma(p+1)}{N!}
 \sum_{j=0}^\infty\frac{(N+j-1)!}{j!\,(N+2j)^p},
 \qquad \al=N\gm^{1/(p-1)} .
\end{equation}
Since \(p>N\), the series converges. One of the main results of this paper is the following.
\begin{theorem}\label{th:main}
Let \(N\ge1,\ 0<s<1,\ sp=N\). There is \(C_{N,s}<\infty\) such that for every measurable \(\Omega\subset\R^N\) of finite positive measure,
\begin{equation}\label{eq:main}
 \sup_{\substack{u\in X(\Omega)\\ [u]_{s,p}\le1}}
 \frac1{|\Omega|}\int_\Omega \exp\big(\al |u|^{p'}\big)\,dx
 \le C_{N,s}.
\end{equation}
Moreover, the constant \(\al\) is sharp in the sense that if we replace it by any \(\alpha>\al\), then the supremum is infinite on every nonempty open set of finite measure, even on smooth compactly supported test functions.
\end{theorem}

The coefficient \(\al\) is exactly \(\alpha^*_{s,N}\) of \cite[Theorem~1.1]{PR}, since \(s/(N-s)=1/(p-1)\). Theorem~\ref{th:main} thus answers the first two of these questions, and it is stronger than Theorem B since it applies to every measurable set of finite measure. In dimension one the series gives
\begin{equation}\label{eq:one}
 \gamma_{s,1}=8\Gamma(p+1)(1-2^{-p})\zeta(p).
\end{equation}
In particular the coefficient for \(p=2\) is \(2\pi^2\) in convention \eqref{eq:seminorm}. On intervals Theorem \ref{th:main} reaches Iula's threshold and includes the endpoint, so it identifies the optimal coefficient throughout the Slobodeckij scale.

The radial kernel \eqref{eq:kernel} also occurs in \cite[Proposition~4.3]{PR}. In logarithmic coordinates the seminorm is an integral of translation differences of a nondecreasing profile \(w\). Its sharp coefficient is determined by the \(p\)-th moment of the translation kernel. Averaging \(w\) on forward intervals gives a majorant with a derivative controlled by the differences that start inside the ball. The remaining differences control the cost of extending the majorant to a zero boundary value. The two estimates use disjoint parts of the energy, which allows the one-dimensional Moser--Jodeit inequality to be applied at the endpoint. More precisely, Theorem~\ref{lem:avg} gives an endpoint criterion for arbitrary nonnegative translation kernels.

\subsection{Singular, whole-space and affine inequalities}
As applications of Theorem \ref{th:main}, we obtain the sharp singular, whole-space, Adachi--Tanaka type and concentration-compactness inequalities. Let \(L=\lceil p-1\rceil\) and define
\begin{equation}\label{eq:trunc}
 \Phi_p(t)=e^t-\sum_{j=0}^{L-1}\frac{t^j}{j!},\qquad t\ge0 .
\end{equation}
For \(\alpha>0\) set
\[
 \mathcal A(\alpha)=
 \sup_{\substack{0\ne u\in W^{s,p}(\R^N)\\ [u]_{s,p}\le1}}
 \frac{\int_{\R^N}\Phi_p(\alpha |u|^{p'})\,dx}{\|u\|_p^p},
\]
and write \(\mathcal C\) for the supremum in \eqref{eq:full} below at \(\alpha=\al\).
\begin{theorem}\label{th:apps}
Let \(N\ge1\), \(0<s<1\) and \(sp=N\).
\begin{enumerate}[label=\textup{(\roman*)},ref=\thetheorem\textup{(\roman*)},leftmargin=*,widest=iii]
\item\label{th:singular} Let \(\Omega\subset\R^N\) be a measurable set of finite positive measure, let \(0\le\beta<N\), let \(B_R=\Omega^\star\) be the ball centered at zero with \(|B_R|=|\Omega|\), and set \(\lambda=N-\beta\). For \(u\in X(\Omega)\) with \([u]_{s,p}\le1\),
\begin{equation}\label{eq:weighted}
 \int_\Omega e^{(1-\beta/N)\al |u|^{p'}}|x|^{-\beta}\,dx
 \le C_{N,s,\beta}\int_{B_R}|x|^{-\beta}\,dx .
\end{equation}
On every open set of finite measure containing zero the coefficient is best possible.
\item\label{th:full} For \(\alpha>0\),
\begin{equation}\label{eq:full}
 \sup_{\substack{u\in W^{s,p}(\R^N)\\ [u]_{s,p}^p+\|u\|_p^p\le1}}
 \int_{\R^N}\Phi_p(\alpha |u|^{p'})\,dx <\infty
 \quad\Longleftrightarrow\quad \alpha\le\al .
\end{equation}
\item\label{th:at} There holds
\begin{equation}\label{eq:equiv}
 \mathcal C=\sup_{0<b<1}\frac{1-b}{b}
       \mathcal A\!\left(\al b^{1/(p-1)}\right).
\end{equation}
The supremum \(\mathcal A(\alpha)\) is finite exactly for \(0<\alpha<\al\). As \(\alpha\uparrow\al\),
\begin{equation}\label{eq:atgrowth}
 \mathcal A(\alpha)\asymp_{N,s}
 \frac{1}{1-(\alpha/\al)^{p-1}} .
\end{equation}
Each finite subcritical supremum has a nonnegative, radially decreasing extremal.
\item\label{th:lions} Let \(0<|\Omega|<\infty\). Suppose \(u_j\in X(\Omega)\),
\([u_j]_{s,p}\le1\), and \(u_j\rightharpoonup u\) in \(W^{s,p}(\R^N)\). Put
\(\Lambda(u)=\al(1-[u]_{s,p}^p)^{-1/(p-1)}\), with value \(+\infty\) if \([u]_{s,p}=1\). For \(0<\alpha<\Lambda(u)\),
\begin{equation}\label{eq:lions}
 \sup_j\int_\Omega e^{\alpha |u_j|^{p'}}\,dx<\infty .
\end{equation}
For every \(a\in(0,1)\) and every nonempty open \(\Omega\) of finite measure, there is such a sequence with \([u]_{s,p}=a\) for which the integrals at \(\alpha=\Lambda(u)\) tend to infinity.
\end{enumerate}
\end{theorem}
For the first-order counterpart of Theorem \ref{th:singular} see Adimurthi and Sandeep \cite{AS}. Theorem \ref{th:at} is a version of the equivalence of Lam, Lu and Zhang \cite{LLZ}, obtained from the critical scaling. For $\alpha$ below the optimal exponent $\alpha^*_{\rm opt}$ of Parini and Ruf, the finiteness of $\mathcal A(\alpha)$ and the existence of extremals, which are radially decreasing about some point, were proved by Chen, Lu and Zhang \cite[Theorems~1.6--1.8]{ChenLuZhang21}. By Theorem \ref{th:main}, $\alpha^*_{\rm opt}(s,\Omega)=\al$ for every bounded Lipschitz domain $\Omega$ when $N\ge2$, so Theorem \ref{th:at} gives the exact range of finiteness of $\mathcal A$, together with the growth \eqref{eq:atgrowth} and the identity \eqref{eq:equiv}. For the full norm on $\R^N$ and coefficients below $\alpha^*_{\rm opt}$, a concentration-compactness principle was given by Chen, Lu and Zhang \cite[Theorem~1.9]{ChenLuZhang21}; Theorem \ref{th:lions} is the sharp version on sets of finite measure, in the spirit of Lions \cite{Lions}.

Finally, we consider affine energies. With \(d\theta\) denoting surface area measure, define for \(u\in W^{s,p}(\R^N)\)
\begin{equation}\label{eq:aff}
 d_u(\theta)=\int_0^\infty\|u(\cdot+r\theta)-u\|_p^p\,\frac{dr}{r^{1+N}},
 \qquad
 \En(u)^p=\frac{S_N^2}{\displaystyle\int_{\sphere}d_u(\theta)^{-1}\,d\theta}
\end{equation}
when \(u\ne0\), and put \(\En(0)=0\). This is the symmetric fractional polar projection body energy of Haddad and Ludwig \cite{HaddadLudwig}, normalized to agree with \eqref{eq:seminorm} on radial functions. It is invariant under volume preserving linear changes of variables, and it is bounded by \([u]_{s,p}\); see Section \ref{sec:aff}.
\begin{theorem}\label{cor:aff}
Theorems~\ref{th:main}, \ref{th:singular} and~\ref{th:full} hold with \(\En(u)\) in place of \([u]_{s,p}\), with the same sharp exponents. On every ball the affine and Euclidean suprema of the unweighted functional agree, and on every ball centered at the origin the same holds for the singular functional; the same holds for the two whole-space full-norm suprema.
\end{theorem}

\subsection{Concentration and critical extremals}
Let \(M=\gamma_{s,N}\) and \(b=M^{-1}\int_0^\infty h^{p+1}k_N(h)\,dh\), where \(k_N\) is the radial kernel in \eqref{eq:kernel}. For \(0<\lambda\le N\) write \(\beta=N-\lambda\), and let \(\mu_\lambda\) be the probability measure \(|x|^{-\beta}dx/\int_{B_1}|x|^{-\beta}dx\) on \(B_1\). Set
\[
 \mathcal S_\lambda=\sup\Bigl\{\int_{B_1}\exp\bigl(\lambda M^{1/(p-1)}|u|^{p'}\bigr)d\mu_\lambda:\ u\in X(B_1),\ [u]_{s,p}\le1\Bigr\}.
\]
By rearrangement one can take a nonnegative radial decreasing maximizing sequence.
Let \(C_\lambda^0\) be the supremum of the limsups of the same integrals over
all admissible radial decreasing sequences converging weakly to zero in \(W^{s,p}(\R^N)\).
Write \(H(p)=\exp(\psi(p)+\gamma_{\rm E})\), where \(\psi=\Gamma'/\Gamma\) and \(\gamma_{\rm E}\) is Euler's constant.
\begin{theorem}\label{th:concentration}
For \(N\ge1\) and \(0<s<1\),
\begin{equation}\label{eq:cl}
1+p\exp\bigl(\lambda b/(p+1)\bigr)\ \le\ C_\lambda^0\
\le\ 1+e^{\lambda b} H(p).
\end{equation}
For each \(N\ge2\) there is \(s_N<1\) such that \(\mathcal S_N\) is attained for \(s_N<s<1\).
For every fixed \(N\ge1\) and \(0<s<1\), \(\mathcal S_\lambda\) is attained for all sufficiently small \(\lambda>0\).
The latter assertion gives critical maximizers for the singular inequality on centered balls when \(\beta\) is sufficiently close to \(N\).
\end{theorem}
In Proposition \ref{pr:improved} we will improve both bounds in \eqref{eq:cl}, and in Corollary \ref{cor:inscribed} we will obtain attainment on bounded domains which contain a large ball. The exact value of \(C_\lambda^0\) in the Hilbert case is given in Section \ref{ss:hilbert}.

\subsection{The regional seminorm}
For a domain \(\Omega\) let
\begin{equation}\label{eq:regional-energy}
 R_\Omega(u)=\int_\Omega\int_\Omega\frac{|u(x)-u(y)|^p}{|x-y|^{2N}}\,dx\,dy ,
\end{equation}
and write \(W^{s,p}_0(\Omega)\) for the closure of \(C_c^\infty(\Omega)\) in the intrinsic norm \((\|u\|_{L^p(\Omega)}^p+R_\Omega(u))^{1/p}\). The boundary condition is essential since constants have zero regional energy. We write \(R_B\) when \(\Omega=B\) is a ball.
\begin{theorem}\label{th:regall}
Let \(0<s<1\) and \(sp=N\).
\begin{enumerate}[label=\textup{(\roman*)},ref=\thetheorem\textup{(\roman*)},leftmargin=*,widest=iii]
\item\label{th:regional} Let \(B\) be a ball. If \(N\ge2\), then
\[
 \sup_{\substack{u\in W^{s,p}_0(B)\\ R_B(u)\le1}}
 |B|^{-1}\int_B e^{\alpha |u|^{p'}}\,dx <\infty
 \quad\Longleftrightarrow\quad \alpha\le\al .
\]
The bound depends only on \(N,s\) for \(0\le\alpha\le\al\).
If \(N=1\), the supremum is infinite for every \(\alpha>0\), even on \(C_c^\infty(B)\).
\item\label{th:lipschitz} Let \(N\ge2\) and let \(\Omega\subset\R^N\) be a bounded Lipschitz domain. Then
\[
 \sup_{\substack{u\in W^{s,p}_0(\Omega)\\ R_\Omega(u)\le1}}
 \frac1{|\Omega|}\int_\Omega e^{\alpha |u|^{p'}}\,dx <\infty
 \quad\Longleftrightarrow\quad \alpha\le\al .
\]
\end{enumerate}
\end{theorem}
Thus on bounded Lipschitz domains the regional and the zero-extension critical coefficients coincide in dimensions \(N\ge2\). On general domains Euclidean balls are not available, and we use instead the intrinsic fractional Hardy inequality; Theorem \ref{th:p3-domain} holds on every bounded domain satisfying this inequality.

\subsection{A fixed anisotropy}\label{ss:aniso}
Let $N\ge2$, $0<s<1$, $p=N/s$ and $p'=q=p/(p-1)$. For an origin-symmetric star body $K\subset\R^N$ with continuous positive radial function $\rho_K$, write $\|h\|_K=|h|/\rho_K(h/|h|)$, $J(h)=\|h\|_K^{-2N}$, $Du(y,z)=u(y)-u(z)$ and
\[
 I_K(u)=\iint |Du(y,z)|^p J(y-z)\,dy\,dz.
\]
For $K=B_1$ this is the Gagliardo energy \eqref{eq:seminorm}, and by Theorem \ref{th:main} the sharp coefficient under $[u]_{s,p}\le1$ is $\alpha_{s,N}$. By the affine P\'olya--Szeg\"o principle of Haddad and Ludwig \cite{HaddadLudwig}, the coefficient $\alpha_{s,N}(|K|/\omega_N)^{2/(p-1)}$ works under $I_K(u)\le1$ uniformly over all bodies of a given volume, and it is sharp on ellipsoids (Theorem \ref{cor:aff}). Sharp $L^2$ Caffarelli--Kohn--Nirenberg inequalities associated with the Minkowski functional $\|\cdot\|_K$ were obtained by the second author \cite{Wei26}, and a two-power affine Hardy--Littlewood--Sobolev inequality, with sharp asymptotics and stability, was studied by the second author and Wu \cite{WeiWu26}. Radial test functions give the upper bound $\alpha_{s,N}(|\mathbb S^{N-1}|^{-1}\int\rho_K^{2N})^{1/(p-1)}$ for the best coefficient. So it is natural to ask for the best coefficient for a fixed body $K$, for instance for the cube.

In this paper we will answer this question. On $S=\mathbb S^{N-1}$ we use ordinary surface measure $d\sigma$. For a star body $L$ with positive Lipschitz radial function put
\begin{equation}\label{eq:MKL}
 M(K,L)=\int_S\int_{\R^N}\left|\log\frac{\|y\|_L}{\|\sigma\|_L}\right|^pJ(\sigma-y)\,dy\,d\sigma ,\qquad M_*=\inf_LM(K,L).
\end{equation}
For a bounded Lipschitz domain $\Omega$ let
\[
 R_{\Omega,K}(u)=\iint_{\Omega^2}|u(y)-u(z)|^pJ(y-z)\,dy\,dz,
\]
and let $W^{s,p}_0(\Omega)$ be the closure of $C_c^\infty(\Omega)$ in the norm $(\|u\|_{L^p(\Omega)}^p+R_{\Omega,K}(u))^{1/p}$.
\begin{theorem}\label{th:aniso}
Let $N\ge2$, $0<s<1$, $p=N/s$, let $K$ be an origin-symmetric star body with continuous positive radial function, and let $\alpha_K=NM_*^{1/(p-1)}$.
\begin{enumerate}[label=\textup{(\roman*)},ref=\thetheorem\textup{(\roman*)},leftmargin=*,widest=iii]
\item\label{th:fixed} There exists a constant $C_{N,s,K}$ such that for every measurable $\Omega$ of finite positive measure and every $u\in X(\Omega)$ with $I_K(u)\le1$, there holds
\[
 \int_\Omega e^{\alpha_K|u|^q}\,dx\le C_{N,s,K}|\Omega| .
\]
Moreover, the constant $\alpha_K$ is sharp in the sense that if we replace it by any $\alpha>\alpha_K$, then the supremum is infinite on every nonempty open set of finite measure. The infimum $M_*$ is attained in the Sobolev completion by a logarithmic profile which is unique up to additive constants.
\item\label{th:ellipsoid} There holds
\[
 \alpha_K\ge\alpha_{s,N}\left(\frac{|K|}{\omega_N}\right)^{2/(p-1)},
\]
with equality if and only if $K$ is an ellipsoid.
\item\label{th:regionalK} For every bounded Lipschitz domain $\Omega$ there exists $C$ such that
\[
 \int_\Omega e^{\alpha_K|u|^{p/(p-1)}}\,dx\le C|\Omega|\qquad\text{for all }u\in W^{s,p}_0(\Omega)\text{ with }R_{\Omega,K}(u)\le1 .
\]
Moreover, the constant $\alpha_K$ is sharp.
\end{enumerate}
\end{theorem}
Since $M(K,cL)=M(K,L)$ for every $c>0$, the minimizer is unique only up to dilation, and we normalize it by $\int_S\log\rho_L\,d\sigma=0$. For $K=B_1$, rotational invariance and uniqueness show that $L$ is a centered ball, and the normalization gives $L=B_1$ and $M_*=\gamma_{s,N}$, so Theorem \ref{th:fixed} contains Theorem \ref{th:main} in dimensions $N\ge2$. The proof of Theorem \ref{th:fixed} gives the endpoint inequality with the coefficient $NM(B_1,B_1)^{1/(p-1)}$ without symmetrization; the identification $M(B_1,B_1)=\gamma_{s,N}$ uses the radial computation of Section \ref{sec:sharp}. Instead of symmetrization, the minimizer of \eqref{eq:MKL} gives a function $v$ which satisfies
\[
 \iint |Dv|^{p-2}Dv\,D\psi\,J(y-z)\,dy\,dz=M_*\psi(0)
\]
for all test functions $\psi$, in the spirit of the work of Shafrir and Spector \cite{SS} on higher order embeddings. Pairing this identity with $u$ and applying H\"older's inequality bounds $|u|$ pointwise by a potential of its local energy density, up to a remainder controlled by the energy away from the diagonal, and the theorem of Fontana and Morpurgo \cite{FM} applies to the resulting homogeneous kernel. The localization in Section \ref{sec:local} loses nothing in the exponent.

The proof of Theorem \ref{th:ellipsoid} uses the equality cases of the Riesz rearrangement inequality due to Burchard \cite{Burchard}. In particular the cube has a strict gap for every $0<s<1$. For $K=B_1$, Theorem \ref{th:regionalK} gives the endpoint of the regional inequality on Lipschitz domains; Theorem \ref{th:p3-domain} gives a different proof, based on the intrinsic fractional Hardy inequality. For the square and $s=1/2$ we obtain the following explicit bounds.
\begin{theorem}\label{th:square}
Let $N=2$, $s=1/2$ and $K=[-1,1]^2$. Then
\[
 \left(\frac{16}{\pi^2}\right)^{1/3}<\frac{\alpha_K}{\alpha_{1/2,2}}<1.185162,
\]
while the radial upper bound is $(16/(3\pi))^{1/3}\approx1.192934$.
\end{theorem}
The left endpoint is approximately $1.174735$. The upper bound is proved by an exact evaluation of $M(K,L)$ for a three-mode trial body in terms of $\pi^2$, $\pi$ and $\pi\zeta(3)$.

\subsection{The Hilbert interval}\label{ss:hilbert}
In this subsection $N=1$ and $s=1/2$. Let $I=(-1,1)$ and let $\mathcal H$ be the space of functions $v\in H^{1/2}(\R)$ which vanish outside $I$, with the norm $\|(-\Delta)^{1/4}v\|_{L^2(\R)}$. Iula, Maalaoui and Martinazzi \cite[Theorem~1.1]{IMM} proved the sharp Trudinger--Moser inequality with the constant $\pi$ under the constraint $\|(-\Delta)^{1/4}v\|_{L^2(I)}\le1$. Since $\|(-\Delta)^{1/4}v\|_{L^2(I)}\le\|(-\Delta)^{1/4}v\|_{L^2(\R)}$, there holds
\[
 \sup_{v\in\mathcal H,\ \|(-\Delta)^{1/4}v\|_{L^2(\R)}\le1}\int_Ie^{\pi v^2}\,dx<\infty .
\]
For the Gagliardo seminorm $[u]_{1/2,2}^2=\iint|u(x)-u(y)|^2|x-y|^{-2}dx\,dy$ one has $[u]_{1/2,2}^2=2\pi\|(-\Delta)^{1/4}u\|_{L^2(\R)}^2$, so this bound under the global constraint, together with the sharpness of $\pi$ under the global constraint, is the case $N=1$, $s=1/2$ of Theorem \ref{th:main}; the sharpness follows from the test functions in Section \ref{sec:sharp}. Mancini and Martinazzi \cite{MM} proved by harmonic extension and commutator estimates that the supremum under the global constraint is attained; see Section \ref{sec:attainment}.

The singular version reads as follows. For $0<\lambda\le1$ let $d\mu_\lambda=\frac\lambda2|x|^{\lambda-1}dx$, a probability measure on $I$. By Theorem \ref{th:singular},
\begin{equation}\label{e:Slambda}
 \mathcal S_\lambda=\sup_{v\in\mathcal H,\ \|(-\Delta)^{1/4}v\|_2\le1}\int_Ie^{\pi\lambda v^2}\,d\mu_\lambda<\infty,
\end{equation}
and the constant $\pi\lambda$ is sharp, in the sense that the supremum is infinite for every larger constant. After passing to a subsequence, a maximizing sequence either has a nonzero weak limit, and then it gives a maximizer, or it converges weakly to zero and concentrates. Let $C_\lambda^0$ denote the concentration level, that is, the supremum of the limits of the integrals in \eqref{e:Slambda} along admissible sequences which converge weakly to zero. Since $\lambda\le1$, the weight $|x|^{\lambda-1}$ is symmetric decreasing, so by the Hardy--Littlewood inequality symmetric decreasing rearrangement does not decrease the weighted integrals of the exponential minus one, and it does not increase the energy. Hence we may restrict to even nonincreasing sequences; see Section \ref{sec:hattain} for the weak convergence of the rearranged sequences. Proposition \ref{pr:improved} gives upper and lower bounds for $C_\lambda^0$ in all dimensions. In the Hilbert case we will compute it exactly.

For a bounded, compactly supported probability density $f$ on $\R$, set
\begin{equation}\label{e:B}
 B_\lambda(f)=-\lambda\iint\log|x-y|f(x)f(y)\,dx\,dy-\int f\log f\,dx+(\lambda-1)\int f(x)\log|x|\,dx.
\end{equation}
\begin{theorem}\label{t:hilbert}
Let $0<\lambda\le1$ and $a=\pi\lambda/2$.
\begin{enumerate}[label=\textup{(\roman*)},ref=\thetheorem\textup{(\roman*)},leftmargin=*,widest=iii]
\item\label{t:loghls} There holds
\[
 \sup_fB_\lambda(f)=\log\pi-\frac1a\int_0^a\log\sin\theta\,d\theta ,
\]
where the supremum is taken over bounded compactly supported probability densities. The value is attained in the class of probability densities for which the three integrals in \eqref{e:B} converge absolutely, by
\[
 f_*(x)=\frac{\sin a}{\pi}\,\frac{|x|^{\lambda-1}}{1+2\cos(a)|x|^\lambda+|x|^{2\lambda}} .
\]
\item\label{t:level} The concentration level is
\[
 C_\lambda^0=1+\lambda\,2^{\lambda-1}\exp\Bigl(\sup_fB_\lambda(f)\Bigr)=1+\pi\lambda\,2^{\lambda-1}\exp\Bigl(-\frac1a\int_0^a\log\sin\theta\,d\theta\Bigr) .
\]
In particular $C_1^0=1+2\pi$.
\end{enumerate}
\end{theorem}
Theorem \ref{t:loghls} is a weighted logarithmic HLS type inequality, which may be of independent interest. For $\lambda=1$, Theorem \ref{t:loghls} is the sharp logarithmic HLS inequality on $\R$ of Carlen and Loss \cite{CL}, and $f_*$ is the Cauchy density. For $\lambda<1$ the density $f_*$ is the normalized source of the singular one-dimensional Liouville bubble studied by DelaTorre, Mancini, Pistoia and Provenzano \cite{DMPP}. We will prove Theorem \ref{t:loghls} by a monotone transport argument in the spirit of Calvez and Carrillo \cite{CC2}.

Using the explicit optimizer, we will construct test functions whose integrals exceed the concentration level. Hence we obtain the existence of extremals.
\begin{theorem}\label{t:attain}
For every $0<\lambda\le1$ the supremum $\mathcal S_\lambda$ is attained. More generally, let $I'=(l,r)$ with $l<0<r$. Then the supremum of $\int_{I'}e^{\pi\lambda v^2}|x|^{\lambda-1}dx$ over $v\in H^{1/2}(\R)$ vanishing outside $I'$ with $\|(-\Delta)^{1/4}v\|_2\le1$ is attained. If $0<\lambda<1$, its concentration level for the normalized weight is
\[
 1+\frac{2c^\lambda}{\lambda\int_{I'}|x|^{\lambda-1}dx}\,(C_\lambda^0-1),\qquad c=\frac{-2lr}{r-l}.
\]
If $\lambda=1$, the concentration level for normalized Lebesgue measure is $1+2\pi$ on every bounded interval.
\end{theorem}
For $\lambda=1$ and the centered interval, the existence of an extremal is the theorem of Mancini and Martinazzi \cite{MM}, and Theorem \ref{t:level} identifies its concentration level $1+2\pi$. The value $1+2\pi$ is the analogue of the value $1+e$ of Carleson and Chang \cite{CC} for the disk.

The organization of the paper is as follows. In Section \ref{sec:moser}, we prove a sharp nonlocal version of Moser's one-dimensional lemma. Section \ref{sec:sharp} gives the radial reduction and the proof of Theorem \ref{th:main}. In Section \ref{sec:apps}, we prove Theorems \ref{th:apps} and \ref{cor:aff}, and we study the limit $s\uparrow1$. Section \ref{sec:attainment} establishes Theorem \ref{th:concentration} and the improved bounds for the concentration level, and Section \ref{sec:regional} deals with the regional seminorm and the proof of Theorem \ref{th:regall}. In Section \ref{sec:anisoproofs}, we study a fixed anisotropy and prove Theorems \ref{th:aniso} and \ref{th:square}. Finally, Section \ref{sec:hilbsec} is devoted to the Hilbert interval and the proofs of Theorems \ref{t:hilbert} and \ref{t:attain}.

\section{A sharp nonlocal Moser theorem}\label{sec:moser}
We first isolate the estimate on the logarithmic radius. Let \(p>1\), let \(k:(0,\infty)\to[0,\infty)\) be measurable, and suppose
\begin{equation}\label{eq:moments}
 M=\int_0^\infty h^p k(h)\,dh\in(0,\infty),\qquad
 b=\frac1M\int_0^\infty h^{p+1}k(h)\,dh\in(0,\infty] .
\end{equation}
For a measurable \(w\) put
\begin{equation}\label{eq:J}
 J(w)=\int_{\R}\int_0^\infty |w(t+h)-w(t)|^p k(h)\,dh\,dt .
\end{equation}

\begin{theorem}[Nonlocal Moser theorem]\label{lem:avg}
Let \(M,b\) be as in \eqref{eq:moments}. Fix \(\lambda>0\), and take the supremum of
\(\int_0^\infty\exp(\alpha w(t)^{p'}-\lambda t)\,dt\)
over nonnegative nondecreasing profiles vanishing on \((-\infty,0)\) with \(J(w)\le1\). The supremum is finite for every \(\alpha<\lambda M^{1/(p-1)}\) and infinite for every larger \(\alpha\). At equality it is finite if and only if \(b<\infty\). In that case
\begin{equation}\label{eq:oned}
 \int_0^\infty \exp\big(\lambda M^{1/(p-1)}w(t)^{p'}-\lambda t\big)\,dt
 \le \lambda^{-1} e^{\lambda b} C_p .
\end{equation}
\end{theorem}
\begin{proof}
First suppose \(b<\infty\). We may assume \(w\) bounded, since truncation decreases all differences and the integrals on the left of \eqref{eq:oned} converge monotonically. Let \(J_+\) and \(J_-\) be the contributions to \eqref{eq:J} from \(t>0\) and \(t<0\), respectively. Introduce the probability measure \(d\nu(h)=M^{-1}h^p k(h)\,dh\) and set
\begin{equation}\label{eq:G}
 G(t)=\int_0^\infty\frac1h\int_0^h w(t+y)\,dy\,d\nu(h).
\end{equation}
Monotonicity gives \(G(t)\ge w(t)\). For each \(h>0\), the inner average is locally absolutely continuous with derivative \((w(t+h)-w(t))/h\) almost everywhere. Integrating this identity from \(0\) to \(t>0\) and applying Tonelli to the nonnegative differences gives
\begin{equation}\label{eq:deriv}
 G(t)=G(0)+\int_0^t D(r)\,dr,\qquad
 D(r)=\int_0^\infty \frac{w(r+h)-w(r)}h\,d\nu(h).
\end{equation}
By Jensen's inequality and the definition of \(\nu\),
\begin{equation}\label{eq:inside}
 M\int_0^\infty |D(r)|^p\,dr\le J_+ .
\end{equation}
Thus \(G\) has the required derivative on the positive half-line.

If \(t<0\), the integrand in \eqref{eq:J} vanishes unless \(t+h>0\). Changing variables gives
\begin{equation}\label{eq:outside}
 J_-=\int_0^\infty k(h)\int_0^h w(y)^p\,dy\,dh .
\end{equation}
In the measure \(k(h)\mathbf1_{\{0<y<h\}}\,dy\,dh\), apply H\"older to the product \(w(y)h^{p-1}\). Its integral is \(M G(0)\), and the integral of \(h^{(p-1)p'}\) is \(Mb\). Hence
\begin{equation}\label{eq:cost}
 M\,G(0)^p b^{1-p}\le J_- .
\end{equation}
On \([0,\infty)\) put \(Y=M^{1/p}G\). On \([-b,0]\) interpolate linearly from \(Y(-b)=0\) to \(Y(0)\). By \eqref{eq:inside} and \eqref{eq:cost}, \(\int_{-b}^\infty |Y'|^p\le J_++J_-\le1\).
The Moser--Jodeit inequality \cite{Jodeit}, in the real-exponent formulation of \cite[Theorem~A]{HL}, gives
\(\int_0^\infty \exp(Z(z)^{p'}-z)\,dz\le C_p\)
for nonnegative locally absolutely continuous \(Z\) with \(Z(0)=0\) and \(\int_0^\infty |Z'|^p\le1\). Use \(Z(z)=\lambda^{1/p'}Y(z/\lambda-b)\). Its derivative has the same \(L^p\) integral as \(Y'\). Restricting to \(z\ge\lambda b\), changing variables and using \(G\ge w\) proves \eqref{eq:oned}.

If \(b=\infty\), apply the preceding bound to \(k\mathbf1_{(0,L)}\), whose energy is smaller and whose \(p\)-th moment increases to \(M\). This proves finiteness at each strictly subcritical coefficient. To prove necessity and sharpness put \(w_T(t)=\min(t_+,T)\), \(J_T=J(w_T)>0\). Integrating \(w_T'=\mathbf1_{(0,T)}\) and taking its translation differences in \(L^p\) gives \(J_T\le MT\). On the plateau of \(w_T/J_T^{1/p}\) every supercritical coefficient therefore gives exponential divergence as \(T\to\infty\). For \(h<T\), integration over the plateau and the two linear ends of the increment gives
\[
 \int_{\R}|w_T(t+h)-w_T(t)|^pdt
 =T h^p-\frac{p-1}{p+1}h^{p+1}.
\]
For \(h\ge T\) the same integral is at most \(Th^p\). Thus
\(MT-J_T\ge\frac{p-1}{p+1}\int_0^T h^{p+1}k(h)\,dh\).
At equality of the coefficients the plateau contributes
\[
 \lambda^{-1}\exp\{\lambda T[(MT/J_T)^{1/(p-1)}-1]\}
 \ge\lambda^{-1}\exp\!\left\{\frac{\lambda}{M(p+1)}\int_0^T h^{p+1}k(h)\,dh\right\},
\]
using \((1-x)^{-1/(p-1)}\ge1+x/(p-1)\). This tends to infinity when \(b=\infty\).
\end{proof}

\section{Radial reduction and the sharp coefficient}\label{sec:sharp}
Let \(u^\star\) be the spherical decreasing rearrangement of \(|u|\). The fractional P\'olya--Szeg\"o inequality gives \([u^\star]_{s,p}\le [u]_{s,p}\); see Almgren and Lieb \cite{AL}. Equimeasurability preserves the integral of the exponential minus one. We add back \(|\Omega|\) after placing \(u^\star\) in the centered ball of measure \(|\Omega|\). Dilations preserve \eqref{eq:seminorm} and the normalized integral in \eqref{eq:main}. Thus it suffices to consider a nonnegative decreasing radial \(v\) supported in \(\overline{B_1}\).

Write \(w(t)=v(e^{-t})\), identifying \(v\) with its radial profile. Then \(w\) is nondecreasing and vanishes for \(t<0\). Ordering the two radii in \eqref{eq:seminorm} and substituting \(r=e^{-t}\), \(\rho=e^{-(t+h)}\), we get
\begin{equation}\label{eq:kernel}
 [v]_{s,p}^p=J(w),\qquad
 k_N(h)=2S_N^2\frac{e^{-Nh}(1+e^{-2h})}{(1-e^{-2h})^{N+1}} .
\end{equation}
Indeed, the Jacobian and the kernel homogeneity contribute \(e^{-Nh}\), and the factor two accounts for both orderings. With surface area measure in both variables, the angular integral for \(a=\rho/r<1\) is
\begin{equation}\label{eq:ang}
 \int_{\sphere}\int_{\sphere}|\theta-a\eta|^{-2N}\,d\theta\,d\eta
 = S_N^2\frac{1+a^2}{(1-a^2)^{N+1}} .
\end{equation}
For \(N\ge2\) this is the radial integration formula of \cite[Proposition~4.3]{PR}. It also follows by squaring the ball Poisson kernel with normalized boundary measure: the convolution of the kernels of radii \(a,a\) is the kernel of radius \(a^2\), as follows by harmonic extension of spherical harmonics. Evaluation on the diagonal gives \eqref{eq:ang}. For \(N=1\) sum over the four pairs in \(\{-1,1\}^2\), obtaining \(2(1-a)^{-2}+2(1+a)^{-2}\).

Both moments in \eqref{eq:moments} exist, since \(k_N(h)=O(h^{-N-1})\) at zero, \(p>N\), and \(k_N\) decays exponentially at infinity. To evaluate the first moment, differentiate \(e^{-Nh}(1-e^{-2h})^{-N}\), or expand its binomial series, to obtain
\begin{equation}\label{eq:series}
 k_N(h)=\frac{2S_N^2}{N}\sum_{j=0}^\infty
 \binom{N+j-1}{j}(N+2j)e^{-(N+2j)h}.
\end{equation}
Integration term by term gives \(M=\gm\) in \eqref{eq:constant}. For \(N=1\), summing over odd integers gives \eqref{eq:one}.
Since \(dx=S_N e^{-Nt}dt\) in the radial integral, Theorem~\ref{lem:avg} with \(\lambda=N\) proves \eqref{eq:main}, with \(C_{N,s}=C_p e^{Nb}\).

\begin{proof}[Sharpness in Theorem~\ref{th:main}]
For \(T>0\) set \(w_T(t)=\min\{\max(t,0),T\}\), and let \(v_T\) be its spatial radial profile. Minkowski's inequality gives
\begin{equation}\label{eq:ramp}
 \|w_T(\cdot+h)-w_T\|_{L^p(\R)}\le h\|w_T'\|_p=hT^{1/p},
 \qquad [v_T]_{s,p}^p\le\gm T .
\end{equation}
The normalized \(u_T=v_T/(\gm T)^{1/p}\) satisfies the constraint. On \(B_{e^{-T}}\) it has the constant value \(T^{1/p'}/\gm^{1/p}\). Thus this ball alone contributes at least
\begin{equation}\label{eq:core}
 \omega_N\exp\!\left(NT\left(\alpha/\al-1\right)\right)
\end{equation}
to the integral with coefficient \(\alpha\). This tends to infinity when \(\alpha>\al\). A translation and dilation place the test in a ball compactly contained in the given open set. Convolution with a nonnegative smooth mollifier preserves the seminorm bound by Minkowski and preserves the plateau on a slightly smaller concentric ball, of at least half its former radius if the mollifier radius is sufficiently small for each \(T\). This also proves divergence in the smooth class.
\end{proof}

\begin{corollary}\label{cor:growth}
If \(F:[0,\infty)\to[0,\infty)\) is Borel measurable and \(F(t)\to\infty\), the assertion \eqref{eq:main} fails with \(F(|u|)\exp(\al |u|^{p'})\) on any nonempty open set of finite measure.
\end{corollary}
\begin{proof}
At equality of the coefficients, the contribution in \eqref{eq:core} is \(\omega_N\). The additional factor is constant on the core and tends to infinity. The dilation used in the sharpness proof changes the lower bound only by a fixed positive factor.
\end{proof}
For \(N\ge2\) this answers the endpoint multiplier question in \cite[Section~6]{PR}. In dimension one the obstruction was proved by Iula \cite[Theorem~1.3]{Iula}.

\section{Singular, whole-space and affine inequalities}\label{sec:apps}
In this section we prove Theorems \ref{th:apps} and \ref{cor:aff}, and we study the limit $s\uparrow1$.

\subsection{The singular weight}\label{sec:singular}
\begin{proof}[Proof of Theorem \ref{th:singular}]
Hardy--Littlewood rearrangement applies to \(\mathbf1_{\Omega}(e^{\alpha |u|^{p'}}-1)\) and to \(\mathbf1_{\Omega}\) separately. The fractional rearrangement inequality reduces the energy to the ball. Scaling to \(R=1\), the normalized weighted integral there is
\(\lambda\int_0^\infty \exp(\alpha w(t)^{p'}-\lambda t)\,dt\).
Theorem~\ref{lem:avg} proves \eqref{eq:weighted} for \(\alpha=\lambda\gm^{1/(p-1)}\). A logarithmic test in a centered interior ball has plateau of weighted measure proportional to \(\exp(-\lambda T)\), giving divergence for every larger coefficient. For the affine assertion use \eqref{eq:ps} and the same radial tests.
\end{proof}
This gives the point-singularity factor in the fractional setting of \eqref{eq:seminorm}; for its first-order counterpart see Adimurthi and Sandeep \cite{AS}. The position of the singularity matters for sharpness. On sets a positive distance from zero the bounded weight permits the unweighted exponent.

\subsection{The full norm on Euclidean space}\label{sec:full}
Classical whole-space results include the subcritical inequality of Adachi and Tanaka \cite{AT}, the full-norm inequalities in \cite{Ruf,LiRuf}, and the level-set method of Lam and Lu \cite{LamLu}. The first retained term has degree \(Lp'\ge p\) in \(|u|\). Iula \cite{Iula} and Zhang \cite{Zhang} obtained full-norm estimates in one and higher dimensions with subcritical coefficients, and Chen, Lu and Zhang \cite{ChenLuZhang21} proved the existence of extremals for the subcritical quotient $\mathcal A(\alpha)$ under the Dirichlet constraint $[u]_{s,p}\le1$. The constraint matters here: attainment for the full-norm functional in \eqref{eq:full} is a different question, and for $N=1$, $s=1/2$ the attainability of both suprema on $\R$ was discussed by Takahashi \cite{Takahashi}. The sharp finite-measure inequality gives the endpoint for the full range \(sp=N\).

\begin{proof}[Proof of Theorem \ref{th:full}]
Take \(\alpha=\al\), let \(a=\|u\|_p\), and assume \(u\ne0\). Then \(0<a<1\) because the only \(L^p(\R^N)\) function of zero seminorm is zero. Put \(\Omega_a=\{|u|>a\}\) and \(v=(|u|-a)_+\). Chebyshev's inequality gives \(|\Omega_a|\le1\); contraction gives \([v]_{s,p}^p\le 1-a^p\). On \(\Omega_a\), two-term H\"older gives
\begin{equation}\label{eq:split}
 |u|^{p'}=(v+a)^{p'}\le
 \frac{v^{p'}}{(1-a^p)^{1/(p-1)}}+1,
\end{equation}
by pairing \(((1-a^p)^{1/p},a)\) with \((v/(1-a^p)^{1/p},1)\). Theorem~\ref{th:main} applied to \(v/(1-a^p)^{1/p}\) bounds the integral over \(\Omega_a\) by \(e^{\al}C_{N,s}|\Omega_a|\). On the complement, \(|u|\le a<1\), and the series in \eqref{eq:trunc} bounds the integrand by \(C_{N,s}|u|^p\).

For sharpness use the unnormalized \(v_T\) of \eqref{eq:ramp}. Its \(L^p\) norm satisfies
\begin{equation}\label{eq:lp}
 \|v_T\|_p^p\le S_N\int_0^\infty t^p e^{-Nt}\,dt=:D_{N,p}.
\end{equation}
The functions \(v_T/(\gm T+D_{N,p})^{1/p}\) are admissible. On the plateau, their \(p'\)-th powers equal \(T\gm^{-1/(p-1)}(1+O(T^{-1}))\). For \(\alpha>\al\) the integral over the core therefore diverges, since \(\Phi_p(z)/e^z\to1\). At \(\alpha=\al\) the same computation gives a positive lower bound on the core integral, so an additional factor \(F(|u|)\) as in Corollary~\ref{cor:growth} again forces divergence.
\end{proof}

\subsection{Subcritical suprema and extremals}\label{sec:subcrit}
In the Hilbert case \(N=1,\ s=1/2\), Takahashi \cite{Takahashi} proved subcritical attainment on the line.
\begin{proof}[Proof of Theorem \ref{th:at}]
Fix \(0<b<1\). Given a competitor \(u\) for \(\mathcal A\), set \(v(x)=b^{1/p}u(\lambda x)\) with \(\lambda^N=b\|u\|_p^p/(1-b)\). Then \([v]_{s,p}^p\le b\) and \(\|v\|_p^p=1-b\). Changing variables in the exponential integral gives the lower bound for \(\mathcal C\) in \eqref{eq:equiv}. Conversely, for an admissible \(v\ne0\) in the definition of \(\mathcal C\), take \(b=[v]_{s,p}^p\in(0,1)\) and \(u=b^{-1/p}v\). The quotient defining \(\mathcal A\) bounds the integral for \(v\) by \(b^{-1}(1-b)\mathcal A(\al b^{1/(p-1)})\).

For the matching lower bound use \(u_T=v_T/(\gm T)^{1/p}\) from \eqref{eq:ramp}. Its \(p\)-th power integral is at most \(D_{N,p}/(\gm T)\) by \eqref{eq:lp}.
For \(\al/2\le\alpha\le\al\) and large \(T\), the plateau alone contributes at least
\(c_N\exp(-NT(1-\alpha/\al))\), since \(\Phi_p(z)/e^z\to1\).
Taking \(T=(1-\alpha/\al)^{-1}\) gives \eqref{eq:atgrowth}. At \(\alpha=\al\), taking \(T\to\infty\) proves divergence. Monotonicity gives divergence at all larger coefficients.

Fix \(0<\alpha<\al\) and let \(c=\alpha^{p-1}/\Gamma(p)\) if \(p-1\) is an integer, and \(c=0\) otherwise. Then \(H(t)=\Phi_p(\alpha t^{p'})-ct^p\) is positive for \(t>0\) and \(H(t)=O(t^{p+\delta})\) at zero for some \(\delta>0\).
In particular \(\mathcal A(\alpha)>c\). By rearrangement and spatial dilation choose a maximizing sequence \(u_j\ge0\) radially decreasing with \(\|u_j\|_p=1\). After passage to a subsequence, \(u_j\rightharpoonup u\) in \(W^{s,p}\) and \(u_j\to u\) almost everywhere by local compactness. The radial bound \(u_j(r)\le(\omega_N r^N)^{-1/p}\) gives
\(\int_{|x|>R} H(u_j)\,dx=O(R^{-N\delta/p})\) uniformly for large \(R\). On each ball \(H(u_j)\) is uniformly bounded in \(L^q\) for some \(q>1\): choose \(q\alpha<\beta<\al\) and use \(H(t)^q\le C(1+\Phi_p(\beta t^{p'}))\) and \(\mathcal A(\beta)<\infty\). Vitali's theorem and the tail estimate give
\(\int H(u)=\mathcal A(\alpha)-c>0\).
Since \(\|u\|_p\le1\) and \([u]_{s,p}\le1\), its quotient is
\(c+\|u\|_p^{-p}\int H(u)\ge\mathcal A(\alpha)\). This proves attainment.
\end{proof}

\subsection{A sharp concentration estimate}\label{sec:lions}
\begin{proof}[Proof of Theorem \ref{th:lions}]
We first note that \(e^{c|u|^{p'}}\in L^1(\Omega)\) for each fixed \(u\in X(\Omega)\) and \(c>0\). Indeed \(v_K=(|u|-K)_+\) has seminorm tending to zero by dominated convergence in \eqref{eq:seminorm} and contraction. Apply Theorem~\ref{th:main} to \(v_K/[v_K]_{s,p}\) for large \(K\), using \(|u|^{p'}\le2^{p'-1}(K^{p'}+v_K^{p'})\). The case \(v_K=0\) is immediate.

From any subsequence select a further subsequence with \(u_j\to u\) almost everywhere by local compactness. Apply the Br\'ezis--Lieb lemma \cite{BL} to
\((u_j(x)-u_j(y))/|x-y|^{2N/p}\) on \(\R^{2N}\). This gives
\[
 \limsup_j [u_j-u]_{s,p}^p\le 1-[u]_{s,p}^p.
\]
Choose \(q>1\) and \(\epsilon>0\) so that
\(q\alpha(1+\epsilon)\limsup_j[u_j-u]_{s,p}^{p'}<\al\).
Young's inequality gives
\(|u_j|^{p'}\le(1+\epsilon)|u_j-u|^{p'}+C_\epsilon |u|^{p'}\).
H\"older's inequality with exponents \(q,q'\), Theorem~\ref{th:main} on the residual, and the preceding fixed-function estimate now give a bound along this subsequence. Applying this reasoning to a putative subsequence of divergent integrals proves \eqref{eq:lions}.

For sharpness take a nonnegative smooth \(u\) with \([u]_{s,p}=a\) supported in one of two disjoint balls compactly contained in \(\Omega\). In the second ball insert a translated, dilated copy \(h_T\) of \(v_T/(\gm T)^{1/p}\). Then \([h_T]_{s,p}\le1\) and \(\|h_T\|_p\to0\), hence \(h_T\rightharpoonup0\) in \(W^{s,p}\). Put \(\epsilon_T=T^{-1/2}\), \(b_T=1-(1-\epsilon_T)^p a^p\) for \(T>1\), and \(u_T=(1-\epsilon_T)u+b_T^{1/p}h_T\).
For nonnegative functions of disjoint supports \([f+g]_{s,p}^p\le[f]_{s,p}^p+[g]_{s,p}^p\), by pointwise comparison of differences. Thus \(u_T\) is admissible. At \(\alpha=\Lambda(u)\) the plateau contribution is a fixed positive multiple of
\[
 \exp\!\left(NT\left[\left(\frac{b_T}{1-a^p}\right)^{1/(p-1)}-1\right]\right).
\]
The term in square brackets is bounded below by a positive multiple of \(T^{-1/2}\) for large \(T\), as follows by differentiation at \(\epsilon_T=0\). This proves divergence.
\end{proof}

\begin{corollary}\label{cor:compact}
Let \(\Omega\) be a nonempty open set of finite measure. For \(0<\alpha<\al\), the finite-domain supremum in \eqref{eq:main} with coefficient \(\alpha\) is attained. At \(\alpha=\al\), any maximizing sequence with a nonzero weak limit along a subsequence yields an extremal.
\end{corollary}
\begin{proof}
Theorem~\ref{th:main} bounds the \(L^p\) norms of admissible functions since \(t^p\le C\exp(\al t^{p'})\). Reflexivity and local compactness provide a weakly converging subsequence with almost everywhere convergence. Its weak limit lies in \(X(\Omega)\) and has seminorm at most one. In the subcritical case choose \(q>1\) with \(q\alpha\le\al\); in the critical case with nonzero weak limit choose \(q>1\) with \(q\al<\Lambda(u)\). Theorems~\ref{th:main} and~\ref{th:lions}, respectively, bound the exponentials in \(L^q(\Omega)\). Vitali's theorem proves convergence of the integrals to that of \(u\).
\end{proof}

\subsection{The affine inequality}\label{sec:aff}
We first collect some properties of the affine energy \eqref{eq:aff}. Polar integration gives \(\int_{\sphere}d_u\,d\theta=[u]_{s,p}^p\). For \(u\ne0\), \(d_u\) is finite almost everywhere and has a positive uniform lower bound: approximation by compactly supported functions gives \(\|u(\cdot+z)-u\|_p^p\to2\|u\|_p^p\) as \(|z|\to\infty\). Thus the reciprocal integral in \eqref{eq:aff} is well defined. Extending \(d_u\) to nonzero vectors by the same integral gives \(d_u(z)=|z|^N d_u(z/|z|)\). The volume of \(\{z:d_u(z)\le1\}\) is \(N^{-1}\int_{\sphere}d_u^{-1}\,d\theta\), and this set transforms by \(T^{-1}\) on replacing \(u\) by \(u\circ T\), \(T\in SL(N)\).

\begin{proof}[Proof of Theorem \ref{cor:aff}]
We record the endpoint argument because the standing hypotheses in the projection-body section of \cite{HaddadLudwig} impose \(p<N/s\). For an origin-symmetric star body \(K\) with continuous positive radial function \(\rho_K\), polar integration gives
\begin{equation}\label{eq:IK}
 I_K(u):=\int_{\R^N}\int_{\R^N}\frac{|u(x)-u(y)|^p}{\|x-y\|_K^{2N}}\,dx\,dy
       =\int_{\sphere}d_u(\theta)\rho_K(\theta)^{2N}\,d\theta .
\end{equation}
If \(|K|=\omega_N\), then \(\int\rho_K^N=S_N\). Cauchy--Schwarz yields \(I_K(u)\ge\En(u)^p\), with equality when \(\rho_K^N\) is proportional to \(1/d_u\). To verify that the optimizing star body is admissible, \(q=d_u^{1/p}\) is even, \(s\)-homogeneous and subadditive by Minkowski's inequality for translation differences. Integrating \(d_u(z)\le 2^{p-1}(d_u(w)+d_u(z-w))\) over \(w\in B_1\) for \(|z|\le1\) bounds \(d_u(z)\) by a multiple of \(\int_{B_2}d_u(w)\,dw<\infty\). Thus \(|q(z)-q(w)|\le C|z-w|^s\); the positive lower bound on the sphere was established above. For \(N\ge2\), \cite[Theorem~12]{HaddadLudwig} gives \(I_K(u)\ge I_K(|u|)\ge I_{B_1}(u^\star)\). Taking the infimum in \(K\), and choosing \(K=B_1\) for the reverse energy comparison, gives
\begin{equation}\label{eq:ps}
 [u^\star]_{s,p}=\En(u^\star)\le\En(|u|)\le\En(u)\le [u]_{s,p}.
\end{equation}
When \(N=1\) the two directional energies are equal by translation, so \(\En(u)=[u]_{s,p}\); ordinary rearrangement gives \eqref{eq:ps}.

On finite-measure sets rearrange to the ball of equal volume and apply Theorem~\ref{th:main}. On all of \(\R^N\) rearrangement also preserves the \(L^p\) norm and the integral of \(\Phi_p\), so Theorem~\ref{th:full} applies. On a ball, for the unweighted functional, and on a ball centered at the pole, for the singular functional, rearrangement to the centered ball of the same volume and \eqref{eq:ps} give the reverse comparison of suprema; the Euclidean admissible class is contained in the affine one. The same holds on the whole space. For the singular functional on a ball which is not centered at the pole, this argument only gives $S_{\rm Euc}(B)\le S_{\rm aff}(B)\le S_{\rm Euc}(B^\star)$, where $B^\star$ is the ball of the same radius centered at the pole. Radial tests in \eqref{eq:ramp}--\eqref{eq:lp} give sharpness in both settings.
\end{proof}

The comparison \eqref{eq:ps} also identifies the affine and Euclidean Adachi--Tanaka suprema at each coefficient, since rearrangement preserves their numerators and denominators.

For a body of arbitrary volume, scaling \eqref{eq:IK} to \(|K|=\omega_N\) shows that
\(I_K(u)\ge (|K|/\omega_N)^2\En(u)^p\). Thus under \(I_K(u)\le1\) the coefficient \(\al(|K|/\omega_N)^{2/(p-1)}\) works uniformly over bodies of a given volume, and it is sharp on ellipsoids by a linear change of variables. For a fixed body \(K\) the sharp coefficient is determined in Theorem \ref{th:fixed}; by Theorem \ref{th:ellipsoid} it is strictly larger than \(\al(|K|/\omega_N)^{2/(p-1)}\) unless \(K\) is an ellipsoid.

\subsection{The first-order limit}\label{sec:limit}
In the Bourgain--Brezis--Mironescu formula \cite{BBM}, the normalization for the double integral matters. Set
\begin{equation}\label{eq:KN}
 K_N=\frac1N\int_{\sphere}|\theta_1|^N\,d\theta
       =\frac{2^{1-N}S_N}{N},\qquad
 a_s=\frac{1-s}{K_N}.
\end{equation}
The angular integral follows from the beta integral and the duplication formula for \(\Gamma\); for \(N=1\) it is the sum over the two directions.
\begin{proposition}\label{pr:limit}
Let \(p=N/s\) and let \(s\uparrow1\). For \(u\in C_c^\infty(\R^N)\),
\begin{equation}\label{eq:bbm}
 a_s [u]_{s,p}^p\longrightarrow \int_{\R^N}|\nabla u|^N\,dx.
\end{equation}
With the normalized constraint \(a_s[u]_{s,p}^p\le1\), the sharp exponential coefficient is \(\widetilde\alpha_{s,N}=\al a_s^{1/(p-1)}\). For \(N\ge2\),
\begin{equation}\label{eq:limit}
 \lim_{s\uparrow1}\widetilde\alpha_{s,N}
       =N S_N^{1/(N-1)} .
\end{equation}
For \(N=1\), \(\lim_{s\uparrow1}\widetilde\alpha_{s,1}^{1/p'}=2\).
\end{proposition}
\begin{proof}
We give the varying-exponent argument for \eqref{eq:bbm}. The part of \eqref{eq:seminorm} with \(|x-y|\ge1\) is bounded uniformly as \(p\downarrow N\). In polar coordinates the remaining part equals
\begin{equation}\label{eq:bbmsplit}
 \int_{\sphere}\int_0^1 r^{p-N-1}
   \int_{\R^N}\left|\frac{u(x+r\theta)-u(x)}{r}\right|^p dx\,dr\,d\theta .
\end{equation}
As \((r,p)\to(0,N)\), the inner integral tends to \(\int|\partial_\theta u|^N dx\) uniformly in \(\theta\), by uniform convergence of the difference quotients, boundedness and common compact support. The probability measures \((p-N)r^{p-N-1}dr\) on \((0,1)\) concentrate at zero. Since \(1-s=(p-N)/p\), multiplication of \eqref{eq:bbmsplit} by \(1-s\) and integration in \(\theta\) give \eqref{eq:bbm}.

From \eqref{eq:kernel}, \(h^{N+1}k_N(h)\to 2^{1-N}S_N^2\) as \(h\downarrow0\). The same probability measures on \(0<h<1\), together with uniform exponential control for \(h\ge1\), give
\begin{equation}\label{eq:momentlimit}
 (p-N)\gm\longrightarrow 2^{1-N}S_N^2,\qquad
 a_s\gm\longrightarrow S_N.
\end{equation}
By homogeneity the normalized sharp coefficient is \(N(a_s\gm)^{1/(p-1)}\), which proves \eqref{eq:limit}. At \(N=1\), taking the \(1/p'\)-th power gives \((a_s\gm)^{1/p}\to2\), in agreement with \(2\|u\|_\infty\le\|u'\|_1\) for compactly supported functions on the line.
\end{proof}

\section{Concentration and critical maximizers on balls}\label{sec:attainment}
\begin{proof}[Proof of Theorem \ref{th:concentration}]
We use two facts from the one-dimensional concentration analysis of Carleson and Chang \cite{CC}, in the real-exponent form of
\cite[Sections~1--2]{HL}. With \(Z(0)=0,\ Z'\ge0,\ \int_0^\infty |Z'|^p\le1\), write
\(F_p(Z)=\int_0^\infty\exp(Z^{p'}-z)\,dz\). Then
\begin{equation}\label{eq:CCfacts}
 Z_j\to0\text{ locally uniformly}\ \Longrightarrow\ \limsup_j F_p(Z_j)\le1+H(p),
 \qquad \sup_Z F_p(Z)>1+H(p).
\end{equation}
To recall precisely the concentration implication, take \(a_j\) to be the first \(z\ge1\) where
\(Z_j(z)^{p'}-z=-2\log z\). Local convergence implies \(a_j\to\infty\) if such points exist,
since \(z-2\log z>0\) on \([1,\infty)\). Before this crossing the integrand is bounded by
\(z^{-2}\) on \([1,a_j)\); the integral up to \(a_j\) tends to \(1\) by dominated convergence.
If there is no crossing this proves convergence of the full integral to \(1\).
For the remaining tail the estimates of \cite[Section~1, Lemmas~1.2.a--b and proof of part~(c)]{HL} give
\(\limsup_j\int_{a_j}^\infty\exp(Z_j^{p'}-z)\,dz\le H(p)\).
They use only the crossing identity and the bound on the remaining derivative energy
\(\int_{a_j}^\infty |Z_j'|^p\le1-Z_j(a_j)^p/a_j^{p-1}\), which follows by H\"older.
These estimates extend from piecewise smooth to locally absolutely continuous increasing profiles:
interpolate linearly on finer finite meshes including \(0,a_j\), and extend constantly beyond the
last mesh point. Jensen preserves the energy bound and Fatou passes the tail bound to the limit.
The strict inequality in \eqref{eq:CCfacts} is the test-function estimate of \cite[Section~2]{HL}.

Suppose \(u_j\rightharpoonup0\) is an admissible radial decreasing sequence.
Compactness in \(L^p(B_1)\) and radial monotonicity imply \(w_j(t)=u_j(e^{-t})\to0\)
at every fixed \(t>0\), choosing monotone representatives. The construction of the forward
average gives
\[
0\le w_j(t)\le G_j(t)\le M^{-1/p}(t+b)^{1/p'}\qquad(t\ge0).
\]
For unbounded profiles the same conclusion and the derivative energy bound follow by applying the
construction to truncations, then passing to the limit by monotone convergence and local weak lower semicontinuity.
In \eqref{eq:G} the displacement \(y\) is uniform on \([0,h]\) given \(h\), and \(h\) has law \(\nu\); hence its mean is
\(\int_0^\infty(h/2)\,d\nu(h)=b/2<\infty\), and the bound above makes \(w_j(t+y)\) dominated by an integrable function of \(y\).
Thus dominated convergence in the average gives \(G_j(t)\to0\) for each fixed \(t\);
monotonicity makes this uniform on compact intervals. The rescaled extensions
\(Z_j(z)=\lambda^{1/p'}Y_j(z/\lambda-b)\) used in the proof of the one-dimensional inequality
also tend to zero locally uniformly. The constant term of the exponential can be subtracted before making the change
of variables. Since the extended integrand \(e^{-z}(\exp(Z_j^{p'})-1)\) is nonnegative,
\[
 \int_{B_1}\!\bigl(e^{\lambda M^{1/(p-1)}|u_j|^{p'}}-1\bigr)\,d\mu_\lambda
 \le e^{\lambda b}\bigl(F_p(Z_j)-1\bigr).
\]
This proves the upper bound in \eqref{eq:cl}.

For the lower bound take \(w_T(t)=\min(t_+,T)\) and divide by \(J(w_T)^{1/p}\).
If \(m=\min(h,T)\) and \(L=\max(h,T)\), integration of the trapezoidal difference gives
\[
 \int_\R |w_T(t+h)-w_T(t)|^p dt=(L-m)m^p+\frac{2m^{p+1}}{p+1}.
\]
Consequently \(J(w_T)=MT-\frac{p-1}{p+1}Mb+o(1)\), with remainder bounded in absolute
value by \(C_p\int_T^\infty h^{p+1}k_N(h)\,dh\). The normalized spatial profiles
are bounded in \(W^{s,p}\) and tend to zero in \(L^p\), hence converge weakly to zero.
In the integral \(\lambda\int_0^\infty \exp(\lambda M^{1/(p-1)}
w_T(t)^{p'}/J(w_T)^{1/(p-1)}-\lambda t)\,dt\), the portion \(t\ge T\) tends to
\(e^{\lambda b/(p+1)}\).
On \([0,T/2]\) the limit contribution is \(1\); on \([T/2,T]\) set \(t=T-y\).
There the exponent tends to \(\lambda b/(p+1)-\lambda y/(p-1)\), giving contribution
\((p-1)e^{\lambda b/(p+1)}\).
For both passages to the limit, use
\(x-x^{p'}\ge c_p\min(x,1-x)\) for \(0\le x\le1\): the exponent on \([0,T]\)
is bounded by \(O(1)-c_p\lambda\min(t,T-t)\). This proves \eqref{eq:cl}.

For the comparison with the local supremum take any bounded nondecreasing piecewise smooth \(Z\)
in \eqref{eq:CCfacts} having bounded, compactly supported derivative. Set
\(w(t)=M^{-1/p}\lambda^{-1/p'}Z(\lambda t)\) for \(t\ge0\), extended by zero for \(t<0\).
Minkowski gives \(J(w)\le M\int_\R|w'|^p\le1\), and the normalized weighted integral equals \(F_p(Z)\).
Thus \(\mathcal S_\lambda\) is bounded below by the supremum of these local trials, independently of \(\lambda\).
Their supremum equals the one in \eqref{eq:CCfacts}: truncating the derivative in size
and in support and integrating it from zero gives increasing approximants and preserves the energy;
piecewise linear approximation then follows by Jensen and Fatou.
For fixed \(p\), choose a trial with \(F_p(Z)>1+H(p)\). It exceeds the right side of \eqref{eq:cl}
for all sufficiently small \(\lambda>0\). For \(\lambda=N\) and \(N\ge2\), choose a trial \(Z\)
at \(p=N\) with a strict gap. The functions \(Z/\|Z'\|_p\) retain a strict gap above
\(1+H(p)\), bounded away from zero, for \(p\) in a right neighborhood of \(N\), by dominated
convergence and continuity of \(H\). As \(p=N/s\downarrow N\) we have \(b\to0\):
in the kernel formula the \((p+1)\)-st moments stay bounded, while the \(p\)-th moments diverge
because \(h^{N+1}k_N(h)\) tends to a positive constant at zero. This gives
\(\mathcal S_\lambda>C_\lambda^0\) in both ranges asserted.

For completeness, strict separation implies existence as follows. A rearranged maximizing sequence is bounded
in \(W^{s,p}(\R^N)\) by the endpoint estimate on the ball, hence has a subsequence
converging weakly and almost everywhere to \(u\in X(B_1)\). Strict separation implies \(u\ne0\).
The Br\'ezis--Lieb lemma applied to the difference quotients with pair measure
\(|x-y|^{-2N}\,dx\,dy\) gives \(\limsup [u_j-u]_{s,p}^p\le1-[u]_{s,p}^p<1\).
Choose \(K\) so large that \(u-u^{(K)}\) has arbitrarily small seminorm, where
\(u^{(K)}=\min(u,K)\); this follows by dominated convergence for the difference quotients.
Then \(\limsup [u_j-u^{(K)}]_{s,p}<1\).
Since \(u^{(K)}\) is bounded, the elementary power inequality
\((a+c)^{p'}\le(1+\varepsilon)a^{p'}+C_{\varepsilon,p}c^{p'}\),
and the weighted endpoint estimate applied to the scaled remainder show that
\(\exp(\lambda M^{1/(p-1)}|u_j|^{p'})\) is uniformly bounded in
\(L^{1+\delta}(d\mu_\lambda)\) for some \(\delta>0\). Vitali convergence gives attainment.
\end{proof}

The two bounds in \eqref{eq:cl} can be improved as follows. For $j\ge0$ put $K_j(y)=\int_y^\infty h^jk_N(h)\,dh$, write $q=p'$, and set
\[
 b_*=\frac1M\int_0^\infty K_{p-1}(y)^qK_0(y)^{-1/(p-1)}\,dy .
\]
\begin{proposition}\label{pr:improved}
For $N\ge1$, $0<s<1$ and $0<\lambda\le N$, we have $0<b_*<b$ and
\begin{equation}\label{eq:improved}
 1+\max\Bigl\{p\,e^{\lambda b/(p+1)},\ H(p)e^{\lambda b/[2(p+1)]}\Bigr\}<C_\lambda^0\le1+H(p)e^{\lambda b_*}.
\end{equation}
In particular,
\[
 \lim_{\lambda\downarrow0}C_\lambda^0=1+H(p),\qquad \lim_{s\uparrow1}C_\lambda^0=1+H(N)\quad(N\ge2,\ 0<\lambda\le N\text{ fixed}).
\]
\end{proposition}
\begin{proof}
\emph{Upper bound.} For the forward average \eqref{eq:G},
\[
 G(t)=M^{-1}\int h^{p-1}k(h)\int_0^h w(t+y)\,dy\,dh,
\]
one has $w\le G$ and $M\int_0^\infty |G'|^p\le J_+$, where
$J_+$ is the energy of differences starting at $t\ge0$. Fubini gives
$MG(0)=\int w(y)K_{p-1}(y)\,dy$ and the complementary energy is
$J_-=\int w(y)^p K_0(y)\,dy$. H\"older's inequality yields
$MG(0)^p b_*^{1-p}\le J_-$.
At each $y>0$, strict Jensen for the probability measure
$k(h)\mathbf 1_{h>y}\,dh/K_0(y)$ gives
$K_{p-1}(y)^q K_0(y)^{1-q}<K_p(y)$.
Integrating gives $b_*<b$. Extend $M^{1/p}G$ linearly from zero
at $-b_*$ to its value at zero. The extended function has derivative
$p$-energy at most one. Along a weakly null radial sequence, $w_j\to0$
at each positive $t$ by compactness in spatial $L^p$ and monotonicity.
The extended functions converge locally uniformly to zero. Indeed the
derivative estimate bounds $w_j(t)\le M^{-1/p}(t+b_*)^{1/q}$.
This bound is integrable in the forward average at each fixed $t$
by the moment condition, and monotonicity makes pointwise convergence of
$G_j$ locally uniform. The one-dimensional Carleson--Chang concentration
bound \eqref{eq:CCfacts}, applied
after the change of variable $z=\lambda(t+b_*)$ to the exponential
minus one, gives the upper bound in \eqref{eq:improved}.

\emph{The effect of the entrance.}
Let $R(\tau)=\min(\tau,0)$, and let $V$ be a nondecreasing function
equal to $\tau$ for $\tau<-a$ and bounded above. Write
$r_T=\min(t_+,T)$, $W_T=(T+V(t-T))_+$, and assume that
$J(W_T)<\infty$ for one $T>a+1$. Define
\[
 d(V)=\lim_{T\to\infty}(J(W_T)-MT),\qquad
 L(V)=e^{-\lambda d(V)/((p-1)M)}
       \lambda\int_{\mathbb R}e^{\lambda(qV(\tau)-\tau)}\,d\tau .
\]
These limits and integrals are finite. Indeed, $Q=V-R$ is bounded,
has finite total variation, and vanishes on $(-\infty,-a)$.
For $T>a$, $W_T=r_T+Q(\cdot-T)$. If $h<T-a$, the difference between
the two inner energies of $W_T$ and $r_T$ equals
$\int(|\Delta_hV|^p-|\Delta_hR|^p)$. This follows by translating the
integrand on the support of $\Delta_hQ(\cdot-T)$, where both arguments
of $r_T$ are positive or zero and $r_T(t)=T+R(t-T)$.
For $h\ge1$ the integral of the absolute difference is at most
$C(1+h)^{p-1}h\,\mathrm{Var}(Q)$, by
$\|\Delta_h Q\|_1\le h\,\mathrm{Var}(Q)$ and
$\|\Delta_h Q\|_\infty\le2\|Q\|_\infty$. At $0<h<1$ the
difference is independent of $T>a+1$, and its absolute integral is
finite after integration against $k$ by the assumed finite energy for
one $T$ and the finite energy of its ramp. Consequently
\[
 d(V)=-\frac{p-1}{p+1}Mb+
 \int k(h)\int_{\mathbb R}
         (|\Delta_hV|^p-|\Delta_hR|^p)\,dt\,dh ,
\]
using $J(r_T)-MT\to-(p-1)Mb/(p+1)$, which was shown in the proof of Theorem \ref{th:concentration}.

Choose any nonzero $\phi\in C_c^\infty((0,\infty))$ with $\phi\ge0$.
For $E(t)=t_+$ set
\[
 D_\epsilon=\int k(h)\int_{\mathbb R}
 \bigl(|\Delta_h(E+\epsilon\phi)|^p-|\Delta_hE|^p\bigr)\,dt\,dh .
\]
The integral is finite, and differentiation at zero gives
\[
 D'_0=p\int_0^\infty\phi(t)\int_t^\infty
        (t^{p-1}-h^{p-1})k(h)\,dh\,dt<0.
\]
Indeed, for small $h$ the integrated absolute difference and its
$\epsilon$-derivative are $O(h^p)$ by Lipschitz continuity and compact
support of the perturbation. For $h\ge1$ they are $O((1+h)^{p-1})$.
For each fixed $h$, shifting the first term in
$\int (\Delta_hE)^{p-1}\Delta_h\phi$ gives the displayed derivative.
Fix $\epsilon>0$ small enough that $D_\epsilon<0$ and
$1+\epsilon\phi'\ge0$. Then
\begin{equation}
 C_\lambda^0-1\ \ge   e^{-\lambda D_\epsilon/((p-1)M)}\sup_V L(V)
   \ >\ \sup_V L(V).
 \label{eq:P13}
\end{equation}
To see this, for large $T$ add $\epsilon\phi(t)$ to $W_T$ and normalize
by its own energy to the power $1/p$. These profiles are nonnegative
and nondecreasing. Their unnormalized energies are
$MT+d(V)+D_\epsilon+o(1)$. For fixed $h$ the entrance perturbation
sees exactly $E$ on all intervals with an endpoint in the support of
$\phi$ once $T$ is large. The preceding bounds give dominated
convergence, since $W_T$ is within $\|Q\|_\infty$ of $r_T$ everywhere,
and is exactly $E$ near the perturbation for small $h$.
The spatial $L^p$ norms of the normalized profiles tend to zero by the
bound $W_T(t)\le t+\|Q\|_\infty$ for $t\ge0$ and exponential decay of
the radial volume measure. Energy boundedness and reflexivity thus give
weak convergence to zero. On a bounded interval in $t$ the integrand
converges to $\lambda e^{-\lambda t}$. On a bounded interval in
$\tau=t-T$ it converges to
\[
 \lambda e^{-\lambda(d(V)+D_\epsilon)/((p-1)M)}
            e^{\lambda(qV(\tau)-\tau)} .
\]
The two intervals are disjoint for large $T$. Fatou's lemma followed
by exhaustion proves the first inequality in \eqref{eq:P13}. The supremum is
positive by $V=R$ and finite by the upper bound. This proves the strict
inequality.

\emph{Explicit trials.}
For $V=R$ one has $L(R)=p\exp(\lambda b/(p+1))$.
For a locally absolutely continuous $V$ as above with
$\int| (W_T)'|^p<\infty$, Jensen on $[t,t+h]$ gives
$|\Delta_hW_T(t)|^p\le h^{p-1}\int_t^{t+h}|W_T'|^p$.
Retain the difference between the two sides just on
$0<h<T-a$, $-h<t<0$, where the right side minus the left side
is $h^{p-1}(t+h)-(t+h)^p$. It follows that
\[
 d(V)\le M\int_{\mathbb R}(|V'|^p-\mathbf 1_{\tau<0})\,d\tau
       -\frac{p-1}{2(p+1)}Mb .
\]
Take $A=\lambda/(p-1)$ and $V_*(\tau)=
\tau-A^{-1}\log(1+e^{A\tau})$. Use the admissible approximations
$V_a(\tau)=\tau$ on $\tau<-a$ and
$V_a(\tau)=V_*(\tau)-V_*(-a)-a$ on $\tau\ge-a$.
They converge uniformly to $V_*$ and their local energy differences
converge to that of $V_*$. Substituting $x=e^{A\tau}$ gives
\[
 \lambda\int_{\mathbb R}e^{\lambda(qV_*-\tau)}\,d\tau=1,
 \qquad
 \int_{\mathbb R}(|V_*'|^p-\mathbf 1_{\tau<0})\,d\tau
   =-\frac{\psi(p)+\gamma_{\rm E}}{A}.
\]
For the second identity one may subtract $(1+x)^{-1}$ inside the
$x$-integral with measure $dx/x$ and use
$\int_0^1(1-y^{p-1})/(1-y)\,dy=\psi(p)+\gamma_{\rm E}$.
Hence $\sup_V L(V)\ge H(p)e^{\lambda b/[2(p+1)]}$.
Equation \eqref{eq:P13} proves the strict lower bound in \eqref{eq:improved}. The first
limit follows immediately. For the second, $p\downarrow N$,
the $(p+1)$-st kernel moments stay bounded, and the $p$-th moments
diverge since $k(h)\asymp h^{-N-1}$ at zero. Thus $b\to0$ and the
bounds converge to the stated value.
\end{proof}

The argument proves critical attainment on Euclidean balls in the stated ranges. For \(\lambda=N\)
the equality of the Euclidean and affine suprema on balls gives affine attainment there as well;
volume-preserving linear changes extend this conclusion to ellipsoids.
The value of \(C_\lambda^0\) for general \(s,N\) and critical attainment on general domains are not known. In the Hilbert case \(N=1\), \(s=1/2\), the value of \(C_\lambda^0\) is computed in Theorem \ref{t:level}, and the singular critical supremum is attained there for every \(0<\lambda\le1\) by Theorem \ref{t:attain}.

\paragraph{The Hilbert case on an interval and on the line.}
For \(N=1,\ s=1/2\), Plancherel gives
\([u]_{1/2,2}^2=2\pi\|(-\Delta)^{1/4}u\|_{L^2(\R)}^2\);
indeed \(\int_\R |e^{ih\xi}-1|^2 h^{-2}dh=2\pi|\xi|\).
There is a difference between the norm printed in \cite[(1.3)]{MM}, which is taken on the
interval \(I\), and the global norm displayed for the extremal in \cite[Theorem~1.1]{MM}.
The interval proof in \cite[Section~2.1, (2.1), Propositions~2.10--2.11]{MM}
uses the global norm: subcritical maximizers of the global form satisfy their equation (2.1),
the blow-up alternative bounds the unnormalized excess integral by \(4\pi\) on \((-1,1)\),
and Proposition~2.11 constructs a globally admissible test exceeding \(4\pi\).
Thus their proof gives attainment for the present interval problem by taking \(v=\sqrt{2\pi}\,u\)
and rescaling the interval. On the whole line \(v(x)=\sqrt{2\pi}\,u(2\pi x)\)
identifies our full constraint with \(\|(-\Delta)^{1/4}v\|_2^2+\|v\|_2^2\le1\).
The functionals differ by a fixed factor, so \cite[Theorem~1.2]{MM} also yields
attainment of the full-norm endpoint on \(\R\).

\paragraph{The closure space on a ball.}
On \(B_1\) the space \(X(B_1)\) equals the closure of \(C_c^\infty(B_1)\) in \(W^{s,p}(\R^N)\).
Indeed \(u_r(x)=u(x/r)\), \(0<r<1\), is supported in \(\overline B_r\) and tends to \(u\)
strongly in \(W^{s,p}(\R^N)\) as \(r\uparrow1\); strong continuity of dilations follows by density of
\(C_c^\infty(\R^N)\) and their uniform norm bounds for \(1/2<r<1\).
Mollifying \(u_r\) at radii smaller than \(1-r\) then gives the asserted approximation.
In particular the ball maximizers just constructed belong to the class appearing in \cite{PR}.

\begin{corollary}\label{cor:inscribed}
Let \(\Omega\subset\R^N\) be a nonempty bounded open set, and put
\(\vartheta_\Omega=\sup\{|B|/|\Omega|:\ B\text{ is a ball contained in }\Omega\}\).
Let \(A(p)=\sup_Z F_p(Z)\) over the one-dimensional admissible class in \eqref{eq:CCfacts}.
If \(\vartheta_\Omega(A(p)-1)>e^{Nb}H(p)\), the unweighted critical supremum in
\(X(\Omega)\) is attained. In particular, for each \(N\ge2\) and each such \(\Omega\) with
\[
 \vartheta_\Omega>\frac{H(N)}{A(N)-1},
\]
attainment holds for all \(s\) sufficiently close to \(1\), depending on \(\Omega,N\).
The threshold on the right is strictly smaller than \(1\).
\end{corollary}
\begin{proof}
If an admissible sequence converges weakly to zero in \(X(\Omega)\), it converges strongly
in \(L^p\) by compactness on a containing ball. Its Schwarz rearrangements are bounded in
\(W^{s,p}\) and converge strongly to zero in \(L^p\), hence weakly to zero in \(W^{s,p}\).
The equal-volume ball comparison and \eqref{eq:cl} thus bound the limsup of the normalized
exponential excess (the integral of the exponential minus one divided by \(|\Omega|\)) by \(e^{Nb}H(p)\).
For every ball \(B\subset\Omega\), inserting the local trials from the proof of
Theorem~\ref{th:concentration}, rescaled and translated into \(B\), gives a lower bound
\((|B|/|\Omega|)(A(p)-1)\) for that supremum. The strict gap therefore forces a nonzero
weak limit of a maximizing subsequence. The same difference-quotient splitting, bounded
truncation and Vitali argument as in Theorem~\ref{th:concentration}, now with Lebesgue
measure on \(\Omega\), yields attainment. For the last assertion take a fixed bounded
local trial with bounded compactly supported derivative at \(p=N\) giving strict
separation after multiplication of its excess by \(\vartheta_\Omega\). Its derivative-normalized
functional is continuous in \(p\), while \(e^{Nb}H(p)\to H(N)\), as proved above.
\end{proof}

\section{The regional seminorm}\label{sec:regional}
\begin{proof}[Proof of Theorem \ref{th:regional}]
By scaling it suffices to consider the unit ball. We begin with \(N\ge2\) and a smooth compactly supported test function. Contraction allows us to take its absolute value. Euclidean decreasing rearrangement can increase the regional energy \cite{LiWang}. Instead equip the ball with the Poincar\'e metric, hyperbolic distance \(d\) and volume \(dV=(2/(1-|x|^2))^Ndx\). The distance formula gives
\begin{equation}\label{eq:regional-conformal}
 |x-y|^{-2N}dx\,dy=H(d(x,y))\,dV(x)dV(y),\qquad
 H(r)=(2\sinh(r/2))^{-2N}.
\end{equation}
Let \(v\) be the hyperbolic decreasing rearrangement of \(|u|\) about the origin. It is also decreasing in the Euclidean radius and has support in some ball \(B_{r_0}\), \(r_0<1\), since the support of \(u\) has finite hyperbolic volume. We claim
\begin{equation}\label{eq:regional-rearr}
 R_B(v)\le R_B(u),\qquad
 \int_B e^{\alpha |u|^{p'}}dx\le\int_B e^{\alpha v^{p'}}dx
 \quad(\alpha\ge0).
\end{equation}
For the first assertion truncate \(H\) to \(H_m=\min(H,m)\). Its integral in one variable is finite and independent of the other: at infinity the hyperbolic area density is \(O(e^{(N-1)r})\) and \(H_m(r)=O(e^{-Nr})\).
For \(a,b\ge0\) put \(F(a,b)=a^p+b^p-|a-b|^p\). Then
\[
 F(a,b)=p(p-1)\int_0^a\int_0^b |t-z|^{p-2}dz\,dt .
\]
The hyperbolic Riesz rearrangement inequality \cite[Theorem~2]{BurchardSchmuck} applied to each pair of superlevel sets shows that the integral of \(F(|u(x)|,|u(y)|)H_m(d(x,y))\) increases under rearrangement. The integrals of the two pure powers are preserved. Subtracting and taking the monotone limit \(m\to\infty\) proves the energy comparison. The second assertion follows by layer cake since \(dx/dV\) is radial decreasing: among sets of fixed hyperbolic volume a centered ball has the largest Euclidean volume. Apply this observation to the superlevel sets of \(e^{\alpha |u|^{p'}}-1\).

Let \(w(t)=v(e^{-t})\) for \(t>0\). The radial computation \eqref{eq:kernel} now gives
\[
 R_B(v)=\int_0^\infty\int_0^\infty |w(t+h)-w(t)|^p k_N(h)\,dh\,dt .
\]
We write \(E_0,E_1\) for the contributions from \(0<t<1\) and \(t\ge1\). The profile \(w\) is nondecreasing and vanishes near \(t=0\). We first prove
\begin{equation}\label{eq:regional-anchor}
 \sup_{0<t<1} w(t)+G(1)\le C_{N,s} E_0^{1/p},
\end{equation}
where \(G\) is the forward average in \eqref{eq:G}, with \(M=\gm\).
Indeed \(k_N(h)\ge c h^{-1-N}\) for \(0<h<1\); hence \(E_0\) bounds a fixed multiple of the \(p\)-th power of the one-dimensional \(W^{s,p}(0,1)\) seminorm. The fractional Morrey inequality on the interval gives
\(\|w-w(0)\|_{L^\infty(0,1)}\le C E_0^{1/p}\), since \(sp=N>1\). One may obtain the homogeneous estimate by subtracting the interval mean and applying the fractional Poincar\'e inequality. Further, for \(0<t<1/2\) and \(z\ge2\), \(k_N(z-t)\ge c e^{-Nz}\). Integrating \(w(z)^p\le 2^{p-1}(|w(z)-w(t)|^p+w(t)^p)\) therefore gives
\[
 \int_2^\infty w(z)^p e^{-Nz}dz\le C E_0 .
\]
Monotonicity gives \(w(2)\le C E_0^{1/p}\). By Tonelli,
\(G(1)=\int_0^\infty w(1+y)g(y)dy\), with
\(g(y)=M^{-1}\int_y^\infty h^{p-1}k_N(h)dh\).
This density has unit mass and is bounded by \(C(1+y^{p-1})e^{-Ny}\) for \(y\ge1\). Bound \(w(1+y)\) by \(w(2)\) for \(y<1\). For \(y\ge1\), apply H\"older with the measure \(e^{-Ny}dy\) and the preceding weighted \(L^p\) estimate. This proves \eqref{eq:regional-anchor}.

The derivative argument in \eqref{eq:deriv}--\eqref{eq:inside}, starting at time one, gives \(M\int_1^\infty |G'|^p\le E_1\).
If \(G(1)\le C E_0^{1/p}\), attach to \(M^{1/p}G\) a linear segment on \([1-\ell,1]\) starting at zero, with \(\ell=(MC^p)^{1/(p-1)}\). Its derivative energy is at most \(E_0\). Since \(E_0+E_1\le1\), the Moser--Jodeit inequality after translation and scaling by \(N\), exactly as in \eqref{eq:oned}, gives
\(\int_1^\infty\exp(\al w(t)^{p'}-Nt)dt\le C_{N,s}\).
Equation \eqref{eq:regional-anchor} bounds the integral over \(0<t<1\). Together with \eqref{eq:regional-rearr} this proves the endpoint estimate. For \(u\) in the stated closure, take an approximating smooth sequence \(u_j\) and divide \(u_j\) by \(\max(1,R_B(u_j)^{1/p})\); an almost everywhere convergent subsequence and Fatou's lemma give the bound. Sharpness follows from the compactly supported tests of Theorem~\ref{th:main}, since regional energy is at most the zero-extension energy.

For \(N=1\) take \(B=(0,1)\), \(0<R<1/4\), and
\[
 b_{a,T}(x)=T^{-1}\min\{T,(\log(R/|x-a|))_+\},\qquad a\in\{0,1\}.
\]
The ramp bound \eqref{eq:ramp} gives \([b_{a,T}]_{s,p}^p\le\gm T^{1-p}\) on \(\R\). Mollify each bump with a smooth nonnegative unit-mass mollifier of radius less than \(Re^{-T}/2\). This preserves the bound and the value one near its center. For all large \(T\), the restriction of one minus these two smoothed bumps to \((0,1)\) belongs to \(C_c^\infty(0,1)\), equals one on a fixed interior interval, and has regional seminorm at most \(2\gm^{1/p}T^{-1/p'}\). Normalizing by this upper bound proves divergence.
\end{proof}
We now treat general domains. Let $N\ge2$, $0<s<1$, $p=N/s$ and $p'=p/(p-1)$. Write
\[
 R_D(f)=\int_D\int_D\frac{|f(y)-f(z)|^p}{|y-z|^{2N}}\,dy\,dz,
 \qquad \delta(x)=\operatorname{dist}(x,\mathbb R^N\setminus\Omega).
\]
The space $W^{s,p}_0(D)$ denotes the closure of $C_c^\infty(D)$
in the intrinsic norm $(\|f\|_{L^p(D)}^p+R_D(f))^{1/p}$.
By Theorem \ref{th:regional} and scaling, for every Euclidean ball $B$,
\begin{equation}\label{eq:p3-ball}
 R_B(f)\le 1,\quad f\in W^{s,p}_0(B)
 \quad\Longrightarrow\quad
 \int_B e^{\alpha_{s,N}|f|^{p'}}\,dx\le C_{N,s}|B|.
\end{equation}

\begin{theorem}\label{th:p3-domain}
Let $\Omega\subset\mathbb R^N$ be a nonempty bounded open set
satisfying the intrinsic fractional Hardy inequality
\begin{equation}\label{eq:p3-hardy}
 \int_\Omega |f(x)|^p\delta(x)^{-N}\,dx\le H R_\Omega(f)
 \qquad (f\text{ compactly supported and Lipschitz in }\Omega).
\end{equation}
Then, for $u\in W^{s,p}_0(\Omega)$ and $R_\Omega(u)\le1$, there holds
\[
 \int_\Omega e^{\alpha_{s,N}|u|^{p'}}\,dx\le C_{N,s,H}|\Omega|.
\]
Every larger coefficient has infinite supremum, already on
$C_c^\infty(\Omega)$. In particular the conclusion holds on every
bounded Lipschitz domain.
\end{theorem}
Theorem \ref{th:lipschitz} is the last assertion of Theorem \ref{th:p3-domain}.

We first estimate the energy needed for a continuous function to reach
the same height on every sphere in an annulus. Suppose $h$ is Lipschitz
on $B_{3r}(x)$ and $h=a\ge0$ on a compact subset $K$ of the closed
annulus $r\le |y-x|\le2r$ meeting every concentric sphere therein.
Then
\begin{equation}\label{eq:p3-annulus}
 a^p\le C_{N,p}\left(r^{-N}\int_{B_{3r}(x)}|h|^p
                    +R_{B_{3r}(x)}(h)\right).
\end{equation}
Indeed, critical scaling reduces this to $x=0$, $r=1$.
There is a probability measure $\mu$ supported on $K$ whose pushforward
under $z\mapsto |z|$ is Lebesgue measure on $[1,2]$.
To construct it, pick a point of $K$ at each of $m$ equally spaced
radii and take a weakly convergent subsequence of the empirical
probabilities as $m\to\infty$. In particular
$\mu(B_t(y))\le 2t$ for every $y$. For $z\in K$ and $0<t\le1$, set
$A_t h(z)=|B_t|^{-1}\int_{B_t(z)}h$. Jensen's inequality gives
\[
 |A_{t/2}h(z)-A_t h(z)|^p
 \le C t^{-2N}\int_{B_t(z)}\int_{B_t(z)}|h(y)-h(w)|^p\,dy\,dw
 \le C R_{B_t(z)}(h).
\]
Fubini's theorem and the bound on $\mu(B_t(y))$ yield
$\int_K R_{B_t(z)}(h)\,d\mu(z)\le 2tR_{B_3}(h)$.
Consequently,
\[
 \int_K |A_{t/2}h-A_t h|\,d\mu\le C t^{1/p}R_{B_3}(h)^{1/p}.
\]
The averages converge to $h$ on $K$, and
$\sup_K |A_1h|\le C\|h\|_{L^p(B_3)}$.
Telescoping for $t=2^{-j}$ and summing $\sum_{j\ge0}2^{-j/p}$
proves \eqref{eq:p3-annulus}.

\begin{proof}[Proof of Theorem~\ref{th:p3-domain}]
First take $u\in C_c^\infty(\Omega)$ with $R_\Omega(u)\le1$.
Fix $x\in\Omega$ and write $r=\delta(x)/20$. By continuity the minimum
\[
 a=\min_{r\le t\le2r}\max_{|y-x|=t}|u(y)|
\]
is attained, say at $t=\rho$. Put $h=\min(|u|,a)$,
$v=(|u|-a)_+$ and $\tau=R_\Omega(h)$.
The compact set where $h=a$ in the closed annulus
$r\le |y-x|\le 2r$ meets every concentric sphere there.
On $B_{3r}(x)\subset\Omega$ one has $\delta(y)\le 23r$.
Equations \eqref{eq:p3-annulus} and \eqref{eq:p3-hardy} therefore give
\begin{equation}\label{eq:p3-anchor}
 a^p\le C_{N,p}(1+23^N H)\tau=:D\tau .
\end{equation}
For every pair of points the increments of $h$ and $v$ have the
same sign and sum to the increment of $|u|$. Thus
\[
 R_{B_\rho(x)}(v)\le R_\Omega(v)\le R_\Omega(|u|)-R_\Omega(h)\le1-\tau .
\]
By choice of $\rho$, the Lipschitz function $v$ vanishes on
$\partial B_\rho(x)$ and belongs to $W^{s,p}_0(B_\rho(x))$.
Indeed, the truncations $(v-\varepsilon)_+$ have compact support in that ball
and converge in its intrinsic norm by dominated convergence, after
which each truncation may be mollified.
For $0<\tau<1$, two-term H\"older and \eqref{eq:p3-anchor} give
\[
 |u|^{p'}\le (v+a)^{p'}
 \le \left(\frac{v}{(1-\tau)^{1/p}}\right)^{p'}
       +\left(\frac{a}{\tau^{1/p}}\right)^{p'}
 \le \left(\frac{v}{(1-\tau)^{1/p}}\right)^{p'}+D^{1/(p-1)} .
\]
Applying \eqref{eq:p3-ball} on $B_\rho(x)$ shows
\begin{equation}\label{eq:p3-local}
 \int_{B_r(x)}e^{\alpha_{s,N}|u|^{p'}}\,dy
 \le 2^N C_{N,s}\exp\!\big(\alpha_{s,N}D^{1/(p-1)}\big)|B_r(x)|.
\end{equation}
When $\tau=0$ we have $a=0$ by \eqref{eq:p3-anchor} and use
\eqref{eq:p3-ball} directly. When $\tau=1$, $R_{B_\rho(x)}(v)=0$,
so the zero boundary values imply $v=0$ throughout that ball and
\eqref{eq:p3-anchor} again gives \eqref{eq:p3-local}.

Apply the five-times-radius covering lemma to the family
$\{B_{\delta(x)/100}(x):x\in\Omega\}$. It supplies disjoint balls
$b_j\subset\Omega$ whose fivefold dilates cover $\Omega$. Each $5b_j$
has precisely the radius $r$ used in \eqref{eq:p3-local}.
Their total volume is at most $5^N|\Omega|$.
Sum \eqref{eq:p3-local} to prove the theorem on smooth functions.
For $u\in W^{s,p}_0(\Omega)$ with $R_\Omega(u)\le1$, approximate in the
intrinsic norm by smooth $u_j$, pass to a subsequence converging almost
everywhere, and apply Fatou's lemma to
$u_j/\max\{1,R_\Omega(u_j)^{1/p}\}$.
For a supercritical coefficient use the smooth zero-extension
Moser sequence of Theorem \ref{th:main} in a ball compactly contained in $\Omega$.
Its regional energy is bounded above by its full-space energy.

Finally, Dyda \cite[Theorem~1.1 (T1)]{Dyda}
proves \eqref{eq:p3-hardy} on every bounded Lipschitz domain for
$sp>1$. Here $sp=N\ge2$.
\end{proof}

Thus on bounded Lipschitz domains the regional and the zero-extension critical coefficients coincide in dimensions \(N\ge2\), and the endpoint holds in both settings.

\section{A fixed anisotropy}\label{sec:anisoproofs}
In this section we prove Theorems \ref{th:aniso} and \ref{th:square}.

\subsection{The angular variational problem}\label{sec:angular}
For \(\sigma,\tau\in S\) and \(r\in\mathbb R\) put
\[
 d\mu=e^{-Nr}J(\sigma-e^{-r}\tau)\,dr\,d\sigma\,d\tau,\qquad
 a_\phi=-r+\phi(\sigma)-\phi(\tau).
\]
Then $M(K,L)=\int |a_\phi|^p d\mu$ for \(\phi=\log\rho_L\), since the inner integral in \eqref{eq:MKL} is homogeneous of degree $-N$ in $\sigma$.
For \(|r|\le1\) the density of \(d\mu\) is comparable to
\((r^2+|\sigma-\tau|^2)^{-N}\), since
\(|\sigma-e^{-r}\tau|^2=(1-e^{-r})^2+e^{-r}|\sigma-\tau|^2\).
For \(|r|\ge1\) the density is bounded by \(C e^{-N|r|}\).
Consequently \(r\in L^p(d\mu)\), using \(p>N\) at the origin,
and integration in \(r\) yields a kernel comparable to
\(|\sigma-\tau|^{1-2N}\).
Since \((N-1)+sp=2N-1\), the norm of \(\phi(\sigma)-\phi(\tau)\) in \(L^p(d\mu)\) is equivalent to the \(W^{s,p}(S)\) seminorm. On the mean-zero subspace the fractional Poincar\'e inequality gives coercivity. Weak lower semicontinuity and reflexivity give a minimizer \(\phi\). Strict convexity gives uniqueness up to constants. Since the kernel is invariant under simultaneous reversal of both directions, uniqueness makes the mean-zero minimizer even. Smooth even approximation gives equality of the two infima, also if the trial bodies are required to be origin symmetric. The minimum \(M_*\) is positive because \(a_\phi\) cannot vanish for almost every \(r\) at fixed \(\sigma,\tau\). Also \(sp>N-1\), so the minimizer has a continuous representative. For all \(b\in W^{s,p}(S)\) the Euler equation is
\begin{equation}
 \int |a_\phi|^{p-2}a_\phi\,[b(\sigma)-b(\tau)]\,d\mu=0 .
 \label{eq:P21}
\end{equation}

\subsubsection*{The logarithmic fundamental solution}
Set \(v(e^{-t}\sigma)=t+\phi(\sigma)\).
For every \(\psi\in C_c^\infty(\mathbb R^N)\),
\begin{equation}
 \iint |Dv|^{p-2} Dv\,D\psi\, J(y-z)\,dy\,dz=M_*\,\psi(0).
 \label{eq:P22}
\end{equation}
Indeed, set \(y=e^{-t}\sigma\), \(z=e^{-t-r}\tau\); then \(J(y-z)\,dy\,dz=dt\,d\mu\) and \(Dv=a_\phi\).
Write \(\Psi(t,\sigma)=\psi(e^{-t}\sigma)\). For each pair of directions,
\[
 \int_\mathbb R|\Psi(t,\sigma)-\Psi(t+r,\tau)|\,dt
 \le C_\psi(|\sigma-\tau|+|r|).
\]
This follows by integrating the angular and \(t\) derivatives of \(\Psi\), which have uniformly bounded \(L^1(dt)\) norms. The bound on the right belongs to \(L^p(d\mu)\). H\"older's inequality justifies integration first in \(t\).
Take smooth \(\eta(t)\) equal to zero and one near \(-\infty\) and \(+\infty\), respectively, and let
\[
 b(\sigma)=\int_\mathbb R(\Psi(t,\sigma)-\psi(0)\eta(t))\,dt .
\]
This is smooth in \(\sigma\). The signed \(t\)-integral of the test difference is \(b(\sigma)-b(\tau)-r\psi(0)\).
Equation \eqref{eq:P21} with \(b\), then with \(\phi\), gives
\(\int (-r)|a_\phi|^{p-2}a_\phi\,d\mu=\int|a_\phi|^p d\mu=M_*\), proving \eqref{eq:P22}.

\subsection{Proof of Theorem \ref{th:fixed}}\label{sec:local}
Define
\[
 H(y)=\int_{\mathbb R^N}|v(y)-v(z)|^p J(y-z)\,dz,\qquad
 k(y)=M_*^{-1}H(y)^{1/q}.
\]
Scaling gives \(H(t\sigma)=t^{-N}H(\sigma)\) for \(t>0\), and \(\int_S H(\sigma)d\sigma=M_*\).
Thus
\begin{equation}
 |\{y:k(y)>\lambda\}|=\alpha_K^{-1}\lambda^{-q}
 \qquad(\lambda>0).
 \label{eq:P23}
\end{equation}
For a ball \(D\) and \(x\in D\), write \(T_D f(x)=\int_D k(y-x) f(y)\,dy\). The finite-measure Adams theorem of Fontana and Morpurgo \cite[Theorem 1]{FM}, applied with \(\beta=\beta_0=q\) and \eqref{eq:P23} in both variables, gives
\begin{equation}
 \int_D e^{\alpha_K |T_D g|^q}\,dx\le C |D|,
 \qquad \|g\|_{L^p(D)}\le1 .
 \label{eq:P24}
\end{equation}
Translation and scaling reduce the assertion to one ball.

\subsubsection*{Localization without loss of exponent}
By spatial dilation assume \(|\Omega|=1\). Tile \(\mathbb R^N\) by half-open unit cubes \(Q\). Let \(D,D'\) be concentric balls about the center of \(Q\) with radii \(2\sqrt N,4\sqrt N\). Choose \(\chi\in C_c^\infty(D')\) equal to one on \(D\), with translates of the same cutoff for each cube. Set
\[
 E_D=\iint_{D^2}|Du|^p J,\qquad
 f_D(y)=\left(\int_D |Du(y,z)|^p J(y-z)\,dz\right)^{1/p}
 \quad(y\in D).
\]
Here and below integrals against \(J\) in two variables include \(dy\,dz\). For almost every \(x\in Q\),
\begin{equation}
 |u(x)|\le T_D f_D(x)+ C\big(I_K(u)-E_D\big)^{1/p}.
 \label{eq:P25}
\end{equation}
To prove this, apply the translate of \eqref{eq:P22} at \(x\) to \(\chi u\), initially smooth. On \(D^2\) use H\"older in \(z\) for each \(y\); the \(L^q\)-norm of the flux in that variable is at most \(H(y-x)^{1/q}\). After division by \(M_*\) the result is \(T_D f_D\).
The other relevant pairs form
\[
 \mathcal R=\big[(D'\times\mathbb R^N)\cup(\mathbb R^N\times D')\big]\setminus D^2.
\]
On \(\mathcal R\) at least one distance to \(x\) is bounded above by a fixed \(B\), and one is bounded below by a fixed \(A>0\).
For fixed \(r,\sigma,\tau\) the allowed \(t\)-interval in polar coordinates about \(x\) therefore has length at most \(\log(B/A)+|r|\).
It follows that \(\iint_{\mathcal R}|Dv(\cdot-x)|^p J\le \int(\log(B/A)+|r|)|a_\phi|^p d\mu<\infty\), uniformly in \(x,Q\). The extra moment is finite by the exponential tails and \(\phi\in L^p(S)\).
On the same region the product rule for differences gives
\[
 \iint_{\mathcal R}|D(\chi u)|^p J
   \le C\big(I_K(u)-E_D+\|u\|_p^p\big).
\]
Indeed \(\sup_z\int |\chi(y)-\chi(z)|^p J(y-z)\,dy<\infty\), since \(p>N\).
H\"older on \(\mathcal R\) gives \eqref{eq:P25} with the displayed \(L^p\)-term still included. This intermediate estimate passes to general \(W^{s,p}(\mathbb R^N)\) functions by smooth approximation. In detail, \(f_{D,j}\to f_D\) in \(L^p(D)\) by the reverse triangle inequality, and \(T_D f_{D,j}\to T_D f_D\) in \(L^1(Q)\) since the homogeneous kernel is locally integrable. Take an almost everywhere convergent subsequence.
Finally, for \(u\in X(\Omega)\) choose a fixed radius \(R\) with \(|B_R|>|D|+2\).
For each \(y\), the set \(B(y,R)\setminus(D\cup\Omega)\) has measure at least one, and on it \(J(y-z)\ge c_K R^{-2N}\) and \(u(z)=0\) almost everywhere. These pairs lie outside \(D^2\), giving \(\|u\|_p^p\le C(I_K(u)-E_D)\). This proves \eqref{eq:P25}.

\subsubsection*{Summation over the cubes}
When \(E_D>0\) set \(g=f_D/E_D^{1/p}\), \(t=\alpha_K(T_D g)^q\), and \(b_0=1/(p-1)\). By \eqref{eq:P24}, \(\int_Q e^t\le C\).
For \(I_K(u)\le1\), \eqref{eq:P25} and two-component H\"older give
\[
 |u(x)|\le E_D^{1/p}T_D g+C(1-E_D)^{1/p}
 \le \big((T_D g)^q+C^q\big)^{1/q}.
\]
The translates of \(D\) have bounded overlap, hence \(\sum_Q E_D\le C I_K(u)\). For fixed \(\delta>0\) only \(C/\delta\) cubes have \(E_D\ge\delta\), and the last two bounds control their total contribution.
Choose \(\delta\) small enough that \(2^{q-1}\delta^{b_0}\le1/2\). For \(0<E_D<\delta\), \eqref{eq:P25} implies
\[
 \alpha_K |u(x)|^q \le 2^{q-1}E_D^{b_0}t+C'.
\]
Where \(E_D^{b_0}t\le1\), the exponential is bounded. Otherwise \(e^{\alpha_K|u|^q}\le C e^{t/2}\le C E_D e^t\), since \(e^{-E_D^{-b_0}/2}\le C E_D\).
Thus on each small-energy cube the integral over \(Q\cap\Omega\) is at most \(C(|Q\cap\Omega|+E_D)\). The same conclusion for \(E_D=0\) follows directly from \eqref{eq:P25}. Summation proves the endpoint estimate.

\subsubsection*{Sharpness} Fix a smooth \(\phi\) and let \(L\) have radius \(e^\phi\). Set \(w_T(t)=\min\{\max\{t,0\},T\}\) and \(u_T(x)=w_T(-\log\|x\|_L)\).
In the cylinder variables above, \(I_K(u_T)=\int\int_\mathbb R
 |w_T(t+\phi(\sigma))-w_T(t+r+\phi(\tau))|^p\,dt\,d\mu\).
For a shift \(h\), convexity applied to \(w_T(t+h)-w_T(t)=\int_0^h{\bf1}_{(0,T)}(t+a)\,da\) gives
\(\int|w_T(t+h)-w_T(t)|^pdt\le T |h|^p\).
For fixed \(h\) the quotient by \(T\) tends to \(|h|^p\) by restricting to intervals on which \(w_T\) has slope one. Dominated convergence gives \(T^{-1}I_K(u_T)\to\int |a_\phi|^p d\mu\).
The normalized function takes the value \(T/I_K(u_T)^{1/p}\) on \(e^{-T}L\), of measure \(e^{-NT}|L|\). Therefore the exponential integral tends to infinity if \(\alpha>N(\int|a_\phi|^p d\mu)^{1/(p-1)}\). Translation and dilation fit \(L\) inside any prescribed nonempty open set. Smooth angular approximants to the minimizer give the asserted sharpness for every \(\alpha>\alpha_K\). This completes the proof of Theorem \ref{th:fixed}.

\subsection{Proof of Theorem \ref{th:ellipsoid}}\label{sec:ell}
Write \(B=B_1\), \(\omega_N=|B|\). By Section \ref{sec:sharp}, \(\gamma_{s,N}=M(B,B)\). We prove
\begin{equation}
 \alpha_K\ge \alpha_{s,N}\left(\frac{|K|}{\omega_N}\right)^{2/(p-1)} ,
 \qquad
 \text{with equality precisely for ellipsoids,}
 \label{eq:P26}
\end{equation}
which is Theorem \ref{th:ellipsoid}. Let \(L\) be the star body of the continuous minimizing \(\phi\). The cylinder argument for \(T^{-1}I_K(u_T)\) in the preceding paragraph requires only \(\int|a_\phi|^p\,d\mu<\infty\), so it applies to this \(L\).
For \(a\ge b\ge0\), the integral of \(p(p-1)(t-s)^{p-2}\) over \(b<s<t<a\) is \((a-b)^p\). The superlevel sets of \(u_T\) at heights \(0<t<T\) are the interiors of \(e^{-t}L\); their boundaries have measure zero by polar coordinates. Tonelli and change of scale therefore give
\[
 \frac{I_K(u_T)}{T}=2p(p-1)\int_0^\infty
       (1-h/T)_+ h^{p-2}\,P_K(e^{-h}L,L)\,dh,\qquad
 P_K(E,F)=\int_{E\times F^c} J(y-z)\,dy\,dz .
\]
Passing to the limit gives the same formula for \(M(K,L)\) with the factor \((1-h/T)_+\) omitted.
Write \(E=aL\), \(F=L\), \(a=e^{-h}\), and \(R_H=(|H|/\omega_N)^{1/N}\), \(H^*=R_HB\). Layer cake for \(J\) gives
\[
 P_K(E,F)=\int_0^\infty 2N r^{-2N-1}
 \left(|E|\,|rK|-\int_E\int_F {\bf1}_{rK}(y-z)\,dz\,dy\right)\,dr .
\]
The Riesz rearrangement inequality shows that the expression in parentheses is at least the corresponding expression for \(E^*,F^*,rK^*\). For every \(h>0\), the interval
\[
  (1-a)R_L< r R_K < (1+a)R_L
\]
has positive length and gives strictly admissible triples of radii. Burchard's equality theorem \cite[Theorem~1]{Burchard} implies that equality for such an \(r\) forces \(rK\) to coincide almost everywhere with an ellipsoid. An origin-symmetric star body with continuous radial function agreeing almost everywhere with an ellipsoid is that ellipsoid, by symmetry and continuity. Consequently if \(K\) is not an ellipsoid there is strict inequality after integration on this interval, then after integration in \(h\). The comparison gives
\[
 M_* =M(K,L)>M(K^*,L^*)=R_K^{2N}\gamma_{s,N}.
\]
For \(K=T B\) the substitution \(y=T y'\), \(z=T z'\) reduces the energy to \(|\det T|^2 I_B\), so equality holds in \eqref{eq:P26}. This proves Theorem \ref{th:ellipsoid}.

\subsection{The square}\label{sec:square}
Take \(N=2,\ p=4,\ K=[-1,1]^2\). Write \(w(\theta)=\rho_K(\theta)^4\),
\(\psi(\theta)=\cos(4\theta)\), and \(m(c)=M(K,L_c)\) for the body with radial function \(e^{c\psi}\).
Then
\begin{equation}
 m'(0)=320\pi-96\pi^2>0,\qquad
 \alpha_K<\left(\frac{16}{3\pi}\right)^{1/3}\alpha_{1/2,2}.
 \label{eq:S1}
\end{equation}
In particular a negative fourth harmonic lowers the energy for all sufficiently small amplitudes, and the radial bound is not sharp for the square.

\begin{proof}
The integral defining \(m\) is a polynomial of degree at most four in \(c\), with finite coefficients by H\"older. In cylinder variables, interchange of the two endpoints sends
\((\sigma,\tau,r)\) to \((\tau,\sigma,-r)\), preserves \(d\mu\), and changes the sign of the logarithm. The two linear terms therefore give equal contributions. In the plane this gives
\[
 m'(0)=8\int_0^{2\pi}\!\cos(4\theta)
 \int_{\mathbb R^2}\frac{w(\arg h)}{|h|^4}
       \big(\log|e^{i\theta}+h|\big)^3\,dh\,d\theta .
\]
Both terms before combining them are absolutely integrable: the cube of the logarithm is \(O(|h|^3)\) near \(h=0\), and the logarithmic singularity at \(-e^{i\theta}\) is locally integrable. By rotations and reflection across the real axis the last expression is
\(8J_4\int_0^{2\pi}w(\theta)\cos(4\theta)\,d\theta\), where
\[
 J_4=\int_{\mathbb R^2}\frac{\cos(4\arg h)}{|h|^4}
           \big(\log|1+h|\big)^3\,dh .
\]
The remaining integral is elementary. Set \(y=1+h\) and identify the plane with \(\mathbb C\). The kernel equals \(\operatorname{Re}(1-y)^{-4}\). Its angular average on \(|y|=r\) is one for \(r<1\) and zero for \(r>1\), by the convergent power series in \(y\) and \(y^{-1}\), respectively. Absolute integrability justifies first integrating in the angle, so
\[
 J_4=2\pi\int_0^1 r(\log r)^3\,dr=-3\pi/4.
\]
Finally
\[
 \int_0^{2\pi} w(\theta)\cos(4\theta)\,d\theta
 =8\int_0^1\frac{1-6t^2+t^4}{1+t^2}\,dt
 =16\pi-\frac{160}{3},
\]
by \(t=\tan\theta\) on the first octant. This proves the derivative in \eqref{eq:S1}. Since \(m(0)=\gamma_{1/2,2}\,16/(3\pi)\), the strict exponent bound follows.
\end{proof}

\subsubsection{An exact certificate for a three-mode trial body}
For \(K=[-1,1]^2\), \(s=1/2\), let \(\alpha_0=\alpha_{1/2,2}\).
The trial function
\[
  \phi(\theta)=
  \frac{-769\cos(4\theta)+133\cos(8\theta)-41\cos(12\theta)}{10000}
\]
gives the rigorous numerical bound
\begin{equation}
 \left(\frac{16}{\pi^2}\right)^{1/3}<\frac{\alpha_K}{\alpha_0}
 < 1.185162 .
 \label{eq:S2}
\end{equation}
In fact its cell energy is exactly
\begin{equation}
 M(K,L)=A\pi^2+B\pi+D\pi\zeta(3),
 \qquad
 \begin{aligned}
 A&=\frac{1111908455570317}{193359375000000},\\
 B&=-\frac{210025743724007423467363}{10050665625000000000000},\\
 D&=\frac{13177246}{390625}.
 \end{aligned}
 \label{eq:S3}
\end{equation}
The quotient \((M(K,L)/(6\pi^2\zeta(3)))^{1/3}\) is approximately \(1.1851611134\).

\subsubsection*{Evaluation}
The following finite coefficient formulas give a reproducible calculation of \eqref{eq:S3}. Write \(w(\theta)=\rho_K(\theta)^4\) and
\[
  w_k=\frac{1}{2\pi}\int_0^{2\pi}w(\theta)\cos(k\theta)\,d\theta
      =\frac4\pi\int_0^{\pi/4}\sec^4\theta\cos(k\theta)\,d\theta
      \quad(k=0,4,8,\ldots).
\]
Let \(b_d\) be the cosine coefficients of \(\phi\), and \(S\) the set of their positive frequencies and their negatives. Form the Laurent polynomials
\[
 \left(\sum_{d\in S}\frac{b_{|d|}}2 X^d(1-Z^d)\right)^j
   =\sum_{k,t} C_{j,k,t} X^k Z^t,\qquad 1\le j\le4.
\]
They have rational coefficients. In cylinder variables take \(r>0\), \(\sigma=e^{i\theta},\ \tau=e^{i\theta}Z,\ |Z|=1,\ X=e^{i\theta},\ u=e^{-r}\). The polynomial on the left is the \(j\)-th power of the angular difference. Exchange of the two directions and \(r\mapsto-r\) doubles the positive half of each term of the quartic energy. Averaging in \(\theta\) replaces \(X^kw(\theta+\arg(1-uZ))\) by \(w_{|k|}\exp(-ik\arg(1-uZ))\). Terms with opposite \(k\) are complex conjugates. Consequently
\begin{equation}
 \frac{M(K,L)}{\pi^2}
 =\frac{32\zeta(3)}\pi+
  8\sum_{j=1}^4 \binom{4}{j}(-1)^{4-j}
    \sum_{k\ge0}(2-{\bf1}_{k=0}) w_k B_{j,k}.
 \label{eq:S4}
\end{equation}
Here the constant term is \(M(K,B_1)/\pi^2\), and all the sums in \(k\) are finite.

For \(k>0\), put \(a=k/2-2,\ b=k/2+2\). Multiplying the radial kernel \(|1-uZ|^{-4}\) by the phase just found gives
\((1-u/Z)^a(1-uZ)^{-b}\).
The angular integral is its constant Laurent coefficient after multiplication by \(\sum_t C_{j,k,t}Z^t\). Expanding at \(0<u<1\) and using \(\int_0^\infty e^{-\lambda r}r^{4-j}dr=(4-j)!/\lambda^{5-j}\) yields
\begin{equation}
 B_{j,k}=\sum_t C_{j,k,t}
   \sum_{\substack{0\le i\le a\\ i\ge t}}
   (-1)^i \binom{a}{i}\binom{b+i-t-1}{i-t}
          \frac{(4-j)!}{(2+2i-t)^{5-j}} .
 \label{eq:S5}
\end{equation}
There is no singular limiting operation in \eqref{eq:S5}, since the constant coefficient for each \(k>0\) is a polynomial in \(u\).

For \(k=0\) the constant coefficient for \(Z^t|1-uZ|^{-4}\) is
\[
 u^{|t|}\frac{1+u^2+|t|(1-u^2)}{(1-u^2)^3}.
\]
Set
\[
 P_j(x)=\sum_t C_{j,0,t}x^{|t|/2}\big(1+|t|+(1-|t|)x\big),
 \quad \sum_{\ell\ge0}U_\ell x^\ell=\frac{P_j(x)}{(1-x)^3}.
\]
There is no constant-frequency term for \(j=1\); for \(j=3\) it vanishes after averaging because the Euclidean angular kernel is symmetric in \(\sigma,\tau\), while the cubed difference is antisymmetric. For \(j=2\), \(P_2\) is divisible by \((1-x)^2\). For \(j=4\), \(P_4\) is divisible by \((1-x)^3\). These divisibilities also follow by direct polynomial arithmetic in the definition of \(P_j\) for the profiles evaluated here. Taking \(d=\deg P_j\), termwise radial integration gives
\begin{equation}
 B_{2,0}=\frac14\left(U_d\zeta(3)+
        \sum_{\ell=0}^{d-1}\frac{U_\ell-U_d}{(\ell+1)^3}\right),
 \qquad
 B_{4,0}=\frac12\sum_{\ell=0}^{d-1}\frac{U_\ell}{\ell+1}.
 \label{eq:S6}
\end{equation}
The Fourier coefficients \(w_k\) in \eqref{eq:S4} can be obtained entirely rationally as affine functions of \(1/\pi\). For \(T_k(X)=\sum_m t_m X^m\) the Chebyshev polynomial, use
\[
 w_k=\frac4\pi\sum_m t_m I_{m-4},\quad
 I_{-4}=\tfrac43,\ I_{-2}=1,\ I_0=\tfrac\pi4,\quad
 I_m=\frac{2^{-m/2}}m+\frac{m-1}{m}I_{m-2}\quad(m=2,4,\ldots),
\]
which is integration by parts for \(I_m=\int_0^{\pi/4}\cos^m\theta\,d\theta\).
Substitution into \eqref{eq:S4}--\eqref{eq:S6} gives \eqref{eq:S3}. As an independent check on the degree-four calculation, taking \(\phi=c\cos4\theta\) gives the whole polynomial
\[
\begin{split}
 M(K,L)={}&32\pi\zeta(3)+(320\pi-96\pi^2)c\\
 &+6\pi^2\left[2w_0(4\zeta(3)-3)+\tfrac5{18}w_8\right]c^2\\
 &+4\pi^2\left(\tfrac43w_4+\tfrac1{15}w_{12}\right)c^3
   +\pi^2\left(34w_0-\tfrac45 w_8\right)c^4 .
\end{split}
\]
For example, \(c=-77/1000\) already yields a quotient below \(1.185501\). Finally, \eqref{eq:S2}, and hence Theorem \ref{th:square}, follows from \eqref{eq:P26} and \eqref{eq:S3} using
\(3.1415926535<\pi<3.1415926536\) and
\(1.202056902<\zeta(3)<1.202056905\). For the latter, append to the sum of the first \(1000\) reciprocals of cubes the integral tails \(1/(2\cdot1001^2)\) and \(1/(2\cdot1000^2)\) as lower and upper bounds. The displayed bounds for \(\pi\) also follow from Machin's formula and the alternating arctangent series.

\subsubsection{Numerical minimization}
We also minimized $M(K,L)$ numerically over the bodies with radial function $\rho_L(\vartheta)=\exp(c_4\cos4\vartheta+c_8\cos8\vartheta+c_{12}\cos12\vartheta)$, by Gauss--Legendre quadrature in polar coordinates centered at $\sigma$ and at the origin, with a smooth partition of unity and exponential radial variables. As a check, the quadrature reproduces $M(B_1,B_1)=6\pi^2\zeta(3)=\gamma_{1/2,2}$ and $M(K,B_1)/\gamma_{1/2,2}=16/(3\pi)$ to six digits. Table \ref{tab:sq} lists the results; the row $s=1/2$ is the trial body of \eqref{eq:S3}.
\begin{table}[h]
\caption{The square $K=[-1,1]^2$: values of $M(K,L)/\gamma_{s,2}$.}\label{tab:sq}
\begin{tabular}{@{}cccccc@{}}
\toprule
$s$ & $p$ & $L=B_1$ & best $L$ found & $(c_4,c_8,c_{12})$ & $(|K|/\omega_2)^2$\\
\midrule
$1/2$ & $4$ & 1.6977 & 1.6647 & $(-0.0769,\,0.0133,\,-0.0041)$ & 1.6211\\
$2/3$ & $3$ & 1.6977 & 1.6857 & $(-0.0394,\,0.0042,\,0)$ & 1.6211\\
$4/5$ & $5/2$ & 1.6977 & 1.6940 & $(-0.0192,\,0.0014,\,0)$ & 1.6211\\
\bottomrule
\end{tabular}
\end{table}

In this three-mode family the best bodies found are larger along the diagonals, like $K$ itself, and their gain over the radial value decreases as $s$ increases through the three values in Table \ref{tab:sq}. We have not proved that the minimizer has this shape, or that the optimal gain is monotone in $s$; the upper bound in Theorem \ref{th:square} only uses the trial body. The numerical behavior is consistent with the first-order limit. Indeed, the proof of Proposition \ref{pr:limit}, with the factor $\rho_K(\theta)^{2N}$ kept in the polar integral \eqref{eq:bbmsplit}, gives for $u\in C_c^\infty(\R^N)$
\[
 \lim_{s\uparrow1}a_sI_K(u)=\frac1{NK_N}\int_{\R^N}\int_{\sphere}\rho_K(\theta)^{2N}|\nabla u(x)\cdot\theta|^N\,d\theta\,dx .
\]
For $N=2$ and $K=[-1,1]^2$, square symmetry gives $\int_{\mathbb S^1}\rho_K(\theta)^4\theta\otimes\theta\,d\theta=\frac{16}3\,\mathrm{Id}$, so the limit is $\frac{16}{3\pi}\int_{\R^2}|\nabla u|^2dx$, which is isotropic.

\subsection{The regional inequality with a fixed anisotropy}\label{sec:regK}
In this subsection we prove Theorem \ref{th:regionalK}. Rotational invariance and uniqueness of the mean-zero minimizer give $\phi=0$ for $K=B_1$, hence \(M_*=\gamma_{s,N}\) in this case.

\subsubsection*{A boundary patch estimate}
On the unit cube \(U=(0,1)^N\), if \(w\) vanishes in a neighborhood of the bottom face, then
\begin{equation}
  \int_U |w|^p\le C \iint_{U^2}\frac{|w(y)-w(z)|^p}{|y-z|^{2N}}\,dy\,dz .
 \label{eq:R1}
\end{equation}
We give a direct proof. Partition the bottom face into dyadic cubes \(P\) of side \(h=2^{-j}\) and set \(m_j\) on each \(P\) equal to the average of \(w\) on \(V_P=P\times(0,h)\). Jensen's inequality bounds the \(p\)-th power of the difference of the averages on a child and its parent by \(C h^{-2N}\iint_{V_P^2}|w(y)-w(z)|^p\,dy\,dz\). Since \({\rm diam}\,V_P\le\sqrt N h\), integration in the face variable and summation give
\[
 \|m_{j+1}-m_j\|_{L^p((0,1)^{N-1})}
 \le C 2^{-j(N-1)/p}R_{U,B_1}(w)^{1/p}.
\]
The series converges for \(N>1\), and \(m_j=0\) for all sufficiently large \(j\). This bounds \(|w_U|=|m_0|\) by \(C R_{U,B_1}(w)^{1/p}\). Jensen on \(U\) also controls \(\int_U |w-w_U|^p\) by \(C R_{U,B_1}(w)\), proving \eqref{eq:R1}.

\subsubsection*{Cubes and patches}
Choose a Whitney decomposition of \(\Omega\) into half-open cubes \(Q\) of side \(\ell\) such that \(8\sqrt N\,\ell\le {\rm dist}(Q,\partial\Omega)\le C_N\ell\). This follows, for instance, by taking the maximal dyadic cubes satisfying the lower bound. Write \(D,D'\) for the concentric balls of radii \(2\sqrt N\ell,4\sqrt N\ell\), both in \(\Omega\). For each \(Q\) there is a patch \(A\subset\Omega\setminus D'\) with
\begin{equation}
 |A|\ge c_\Omega\ell^N,\quad {\rm diam}(A\cup D')\le C_\Omega\ell,\qquad
 \ell^{-N}\int_A |u|^p\le C_{\Omega,K} R_{A,K}(u)
 \quad(u\in C_c^\infty(\Omega)).
 \label{eq:R2}
\end{equation}
To construct \(A\), take a boundary point nearest the cube center and use a Lipschitz graph chart there, with uniform slope bound and uniform chart radius \(r_0>0\) by bounded Lipschitz regularity. Above a tangential cube of side \(\rho=c\min(\ell,r_0)\) centered at that point, take the vertical strip between the graph and height \(\rho\) above it. Choose \(c>0\) fixed so the strip lies in the chart and within \(\sqrt N\ell\) of that boundary point. This makes it disjoint from \(D'\). Since \(\ell\le {\rm diam}(\Omega)\), one has \(\rho\ge c_\Omega\ell\), giving the first two bounds. Graph flattening after dilation by \(\rho\) is uniformly bi-Lipschitz and changes the fractional energy by at most a constant, since \(J(h)\) is comparable to \(|h|^{-2N}\). Pull back \eqref{eq:R1}, observing that \(u\) vanishes near the base. The scaling \(sp=N\) gives the final bound.

Put \(E_D=R_{D,K}(u)\) and \(R=R_{\Omega,K}(u)\). Averaging
\(|u(y)|^p\le 2^{p-1}(|u(y)-u(z)|^p+|u(z)|^p)\)
over \(z\in A\), then integrating over \(y\in D'\), gives
\begin{equation}
 \ell^{-N}\int_{D'} |u|^p
 \le C \iint_{D'\times A}|Du|^p J
      + C \ell^{-N}\int_A |u|^p
 \le C (R-E_D).
 \label{eq:R3}
\end{equation}
Both terms on the right of the first inequality use pairs outside \(D^2\).

\subsubsection*{Application of the logarithmic flux}
Take \(\chi\in C_c^\infty(D')\) equal to one on \(D\), chosen by translation and dilation of a single cutoff. For \(x\in Q\), apply \eqref{eq:P22} centered at \(x\) to the zero extension of \(\chi u\). On \(D^2\) this gives the potential \(T_D f_D\) of the preceding section, with \(\|f_D\|_p^p=E_D\). On the remaining pairs \(\mathcal R\) there, the energy of the logarithmic solution is uniformly bounded; dilation changes \(v\) by an additive constant and leaves its energy unchanged. To bound the other factor, on \((D')^2\setminus D^2\) use the product rule for differences and
\(\sup_z\int |\chi(y)-\chi(z)|^p J(y-z)\,dy\le C\ell^{-N}\).
When \(y\in D'\), \(z\notin D'\) the difference is \(\chi(y)u(y)\) and
\[
 |\chi(y)|^p\int_{\mathbb R^N\setminus D'} J(y-z)\,dz\le C\ell^{-N}.
\]
Indeed the integral is at most \(C\,{\rm dist}(y,\partial D')^{-N}\) and
\(|\chi(y)|\le C \min(1,{\rm dist}(y,\partial D')/\ell)\), with \(p>N\).
It follows from \eqref{eq:R3} that
\(\iint_\mathcal R |D(\chi u)|^p J\le C(R-E_D)\). Thus
\[
 |u(x)|\le T_D f_D(x)+C(R-E_D)^{1/p}.
\]
For \(R\le1\) and \(E_D>0\), write \(g=f_D/E_D^{1/p}\). Two-component H\"older yields
\[
  |u(x)|^{p/(p-1)}
   \le (T_Dg(x))^{p/(p-1)}+C^{p/(p-1)} .
\]
For \(E_D=0\), the preceding bound gives \(|u(x)|\le C\). The homogeneous Adams inequality \eqref{eq:P24} therefore yields
\(\int_Q \exp(\alpha_K |u|^{p/(p-1)})\le C |D|\le C'|Q|\).
Sum over the partition. Approximating an arbitrary \(u\in W^{s,p}_0(\Omega)\) in the intrinsic norm by smooth compactly supported functions \(u_j\), divide \(u_j\) by \(\max(1,R_{\Omega,K}(u_j)^{1/p})\) and apply Fatou along an almost everywhere convergent subsequence. This passes the estimate to the prescribed class. For sharpness use the compactly supported ramps of the full-space theorem in an interior ball: their regional energy is bounded by the full energy. They belong to \(W^{s,p}_0(\Omega)\) by interior mollification. This proves Theorem \ref{th:regionalK}.

\section{The Hilbert interval}\label{sec:hilbsec}
In this section we prove Theorems \ref{t:hilbert} and \ref{t:attain}.

\subsection{Preliminaries}\label{sec:hilb}
Throughout this section, $N=1$, $s=1/2$, $p=2$, the letters $J$ and $K_j$ denote the one-dimensional energy and the kernel moments below, and we use the logarithmic radius $t=-\log|x|$. For an even nonincreasing function $u$ on $I$ let $w(t)=u(e^{-t})$, a nondecreasing profile which vanishes for $t<0$. By Section \ref{sec:sharp},
\[
 [u]_{1/2,2}^2=J(w)=\int_0^\infty k(h)\int_{\R}|w(t+h)-w(t)|^2\,dt\,dh,\qquad k(h)=8\frac{e^{-h}(1+e^{-2h})}{(1-e^{-2h})^2}.
\]
Set $M=\int_0^\infty h^2k(h)\,dh=2\pi^2$, $K_j(y)=\int_y^\infty h^jk(h)\,dh$ and
\[
 b_*=\frac1M\int_0^\infty K_1(y)^2K_0(y)^{-1}\,dy,
\]
which is finite by Proposition \ref{pr:improved}. The normalization $v=\sqrt{2\pi}\,u$ gives $\|(-\Delta)^{1/4}v\|_2=1$ if $[u]_{1/2,2}=1$, and $\pi\lambda v^2=\lambda M u^2$.
\begin{lemma}\label{l:forward}
Let $w\ge0$ be nondecreasing, vanishing on $(-\infty,0)$, with $J(w)\le1$. Then there exists a locally absolutely continuous function $Y$ on $[-b_*,\infty)$ with $Y(-b_*)=0$, $\int_{-b_*}^\infty|Y'|^2\le1$ and $M^{1/2}w\le Y$ on $[0,\infty)$. In particular $M^{1/2}w(t)\le(t+b_*)^{1/2}$.
\end{lemma}
\begin{proof}
For the forward average
\[
 G(t)=M^{-1}\int h^{p-1}k(h)\int_0^h w(t+y)\,dy\,dh,
\]
one has $w\le G$ and $M\int_0^\infty |G'|^p\le J_+$, where
$J_+$ is the energy of differences starting at $t\ge0$. Fubini gives
$MG(0)=\int w(y)K_{p-1}(y)\,dy$ and the complementary energy is
$J_-=\int w(y)^p K_0(y)\,dy$. H\"older's inequality yields
$MG(0)^p b_*^{1-p}\le J_-$.
(Here $p=q=2$; the first two assertions are proved in Section \ref{sec:moser}.)
Extend $M^{1/2}G$ linearly from zero at $-b_*$ to its value at zero. The extended function $Y$ has derivative energy at most $J_++J_-\le1$, and $M^{1/2}w\le M^{1/2}G=Y$ on $[0,\infty)$. The last assertion follows by Cauchy--Schwarz.
\end{proof}

\subsubsection*{The Green function}
Let $G$ be the inverse of the Dirichlet half-Laplacian on $I$ with zero exterior values. Its kernel, displayed in \cite[(2.14)]{MM}, is
\[
 \pi G(x,y)=\log\frac{1-xy+\sqrt{(1-x^2)(1-y^2)}}{|x-y|}\le\log\frac2{|x-y|}.
\]
The numerator is at most $2$ by Cauchy--Schwarz. For each bounded probability density $f$ supported in $I$, Hilbert-space duality thus gives, for $\|(-\Delta)^{1/4}v\|_2\le1$,
\begin{equation}\label{eq:P15}
 \pi\left(\int_I vf\right)^2\le \pi\langle f,Gf\rangle
 \le \log2-\iint\log|x-y|f(x)f(y)\,dx\,dy .
\end{equation}

\subsection{The concentration level on the Hilbert interval}\label{sec:hlevel}
In this subsection we prove the first identity in Theorem \ref{t:level}; the second one then follows from Theorem \ref{t:loghls}. We use the normalization $v$ above, with $\|(-\Delta)^{1/4}v\|_2\le1$ and exponent $\pi\lambda v^2$.

\subsubsection*{Upper bound}
Take an admissible weakly null even nonincreasing
sequence $v_j\ge0$. Put $z_j(t)=\sqrt\pi\,v_j(e^{-t})=M^{1/2}w_j(t)$. By Lemma \ref{l:forward}, $z_j\le Y_j$ for an absolutely continuous
$Y_j$ on $[-b_*,\infty)$, with $Y_j(-b_*)=0$ and $\int|Y_j'|^2\le1$.
In particular $z_j(t)\le\sqrt{t+b_*}$. Moreover, $w_j\to0$ locally uniformly on $[0,\infty)$, by compactness in $L^2(I)$ and monotonicity.
Fix $R>2$, sufficiently large that
$t-(2/\lambda)\log t>0$ and increasing for $t\ge R$.
Call $t\ge R$ good when $z_j(t)^2\ge t-(2/\lambda)\log t$.
The bad part of $\int_R^\infty e^{\lambda(z_j^2-t)}dt$ is at most $1/R$.
The good radii eventually leave every compact interval because $z_j\to0$
locally uniformly. We claim that the oscillation of $z_j$ on all good
radii tends to zero.

Indeed let $a\le t$ be two good radii, set $d=t-a$, and write
$c=2/\lambda$, $b=b_*$. Cauchy--Schwarz on $[-b,a]$ gives
$\int_a^\infty|Y_j'|^2\le(b+c\log a)/(a+b)=:\eta$.
Thus
\[
 \sqrt{t-c\log t}\le \sqrt{a+b}+\sqrt{d\eta},
 \qquad
 d(1-\eta)\le b+c\log t+2\sqrt{d(b+c\log a)} .
\]
As $a\to\infty$, one has $\eta\to0$, and $\log t\le\log a+d/a$.
Absorbing the $d/a$ and square root terms yields $d\le C(1+\log a)$,
uniformly. Good values at $a$ and $t$ are therefore both in
$[\sqrt{a-c\log a},\sqrt{a+C(1+\log a)}]$, of length
$O(\log(a)/\sqrt a)$. This proves the claim, independently of whether
the first good radius is attained.

Let $A_j=\{x\in I:-\log|x|\text{ is good}\}$. It suffices to consider
$A_j$ of positive measure. Write
$Z_j=\int_{A_j}e^{\pi\lambda v_j^2}|x|^{\lambda-1}dx$ and let $f_j$
be its normalized integrand extended by zero. The diameter bound for
the good radii places $A_j$ away from zero for each fixed $j$, and the
upper envelope bounds $v_j$ there, so $f_j$ is bounded and compactly
supported. The oscillation claim gives $\mathrm{Var}_{f_j}(v_j)=o(1)$.
The identity
\[
 \log Z_j=\pi\lambda\int v_j^2 f_j-\int f_j\log f_j
                 +(\lambda-1)\int f_j\log|x|
\]
and \eqref{eq:P15} imply $\log Z_j\le\lambda\log2+\sup_f B_\lambda(f)+o(1)$.
On $0\le t<R$ local convergence gives the background contribution in
the limit. Multiplying the good-set integral by $\lambda/2$ for the
probability normalization and sending $R\to\infty$ proves the upper
bound in Theorem \ref{t:level}.

\subsubsection*{Lower bound} Fix any bounded compactly supported probability density $f$,
put $f_\varepsilon(x)=\varepsilon^{-1}f(x/\varepsilon)$ and
$A_\varepsilon=\varepsilon\operatorname{supp} f\subset I$ for small
$\varepsilon$. The kernel formula gives
\[
 \pi m_\varepsilon:=\pi\langle f_\varepsilon,Gf_\varepsilon\rangle
   =\log(1/\varepsilon)+\log2
       -\iint \log|x-y| f(x)f(y)\,dx\,dy+o(1)
\]
by uniform convergence of the numerator on the dilated support.
Set $v_\varepsilon=Gf_\varepsilon/\sqrt{m_\varepsilon}$.
These functions have unit energy and converge weakly to zero: their
energy pairing with a smooth compactly supported test $g$ is
$\int g f_\varepsilon/\sqrt{m_\varepsilon}\to0$, and the tests are
dense in the zero-exterior energy space on $I$. Jensen's inequality
for $f_\varepsilon dx$ and then Cauchy--Schwarz give
\[
 \log\int_{A_\varepsilon}e^{\pi\lambda v_\varepsilon^2}|x|^{\lambda-1}dx
 \ge \pi\lambda\left(\int v_\varepsilon f_\varepsilon\right)^2
    \!-\!\int f_\varepsilon\log f_\varepsilon
       +(\lambda-1)\int f_\varepsilon\log|x|
 =\lambda\log2+B_\lambda(f)+o(1).
\]
Outside $A_\varepsilon$ the exponential is at least one, and the weighted
measure of $A_\varepsilon$ tends to zero. Symmetric decreasing rearrangement
preserves the vanishing $L^2$ norm by compactness, decreases the
energy, and increases the weighted integral by Hardy--Littlewood
(rearrange the exponential minus one). The rearranged sequence is
therefore admissible in the definition of $C_\lambda^0$. This proves
the lower bound, and in particular finiteness of $\sup_f B_\lambda(f)$
follows from the upper bound $C_\lambda^0\le1+H(2)e^{\lambda b_*}$ of Proposition \ref{pr:improved}.

At $\lambda=1$ the sharp logarithmic HLS inequality on $\mathbb R$
of Carlen and Loss \cite{CL} gives $B_1(f)\le\log(2\pi)$, with equality for
$f_0(x)=[\pi(1+x^2)]^{-1}$. The constant also follows by evaluating
at $f_0$: its entropy is $-\log(4\pi)$ and the difference of two
independent standard Cauchy variables has Cauchy scale two, so its
logarithmic energy is $\log2$. The normalized compact truncations of
$f_0$ tend to this value. Entropy and logarithmic energy converge by
dominated convergence, since $f_0$ has finite logarithmic moment and
is bounded, and $|\log|x-y||$ is integrable locally in $x-y$.
This establishes the value $C_1^0=1+2\pi$, which is also the case $\lambda=1$ of Theorem \ref{t:loghls}.

\subsection{The weighted logarithmic inequality}\label{sec:hlog}
In this subsection we prove Theorem \ref{t:loghls}. Let $0<\lambda\le1$ and $a=\pi\lambda/2$, and write $f_*$ for the density in Theorem \ref{t:loghls}.

\subsubsection*{The stationary density}
Putting $t=|x|^\lambda$ shows that $f_*$ has mass one, since
$\int_0^\infty(1+2t\cos a+t^2)^{-1}dt=a/\sin a$.
Its logarithmic potential equals
\begin{equation}\label{eq:weighted-logpot}
 P(x):=\int\log|x-y|\,f_*(y)\,dy
 =\frac1{2\lambda}\log(1+2\cos(a)|x|^\lambda+|x|^{2\lambda}).
\end{equation}
Indeed, on the upper half-plane take
$H(z)=\lambda^{-1}\log(1+(-iz)^\lambda)$ with principal branches.
Then $H'(z)=(-iz)^\lambda/[z(1+(-iz)^\lambda)]$, whose boundary
imaginary part at $x\ne0$ is $-\pi f_*(x)$. Reflect $H'$ to the lower
half-plane by complex conjugation. Applying the Cauchy integral formula
on the two sides of the real axis gives
$H'(z)=\int f_*(y)/(z-y)\,dy$ on the upper half-plane:
the reflected boundary jump is $-2\pi i f_*$, the contribution to the
Cauchy formula on a circle of radius $\epsilon$ is $O(\epsilon^\lambda)$
and the large-circle contribution tends to zero since $H'(z)=O(1/|z|)$.
Hence $\operatorname{Re}H(z)-\int\log|z-y|f_*(y)\,dy$ is constant.
Sending $z=iR$ to infinity gives zero, by the finite logarithmic moment
of $f_*$. Taking boundary values gives \eqref{eq:weighted-logpot} for
$x\ne0$. At $x=0$ both sides are zero by inversion symmetry.
Consequently
\begin{equation}\label{eq:weighted-stationary}
 \log f_*(x)+2\lambda P(x)-(\lambda-1)\log|x|
 =\log(\sin a/\pi)\qquad (x\ne0).
\end{equation}

\subsubsection*{Even reduction and transport}
Symmetric decreasing rearrangement increases $B_\lambda$.
For a compactly supported bounded competitor and its rearrangement,
take $R$ larger than the absolute values of all support points and all
pairwise distances within either support. The kernel
$\log_+(R/|x-y|)$ is symmetric decreasing in $x-y$; Riesz rearrangement
therefore increases its interaction integral. Similarly
$\log_+(R/|x|)$ is symmetric decreasing and its coefficient in the
rewritten weight term is $1-\lambda\ge0$. The entropy is invariant.
We may thus take even $f$. Set $u=\log|x|$, $g(u)=2e^uf(e^u)$ and
$Q(z)=\log(2|\sinh z|)$. Averaging over the two signs gives
\[
 -B_\lambda(f)=-\log2+\mathcal F(g),\qquad
 \mathcal F(g)=\int g\log g+\frac\lambda2\iint Q(u-v)g(u)g(v)\,du\,dv.
\]
Indeed the averaged logarithmic kernel is $(u+v)/2+Q(u-v)/2$,
so the linear moments cancel. Write
$g_*(u)=2e^uf_*(e^u)=\sin(a)/[\pi(\cosh(\lambda u)+\cos a)]$.
The same change of variables in \eqref{eq:weighted-stationary} gives
$\log g_*+\lambda Q*g_*=\mathrm{constant}$.

We first prove the inequality for $g$ of the form
$(1-\epsilon)g_0+\epsilon g_*$, where $0<\epsilon<1$ and $g_0$ is a
smooth compactly supported probability density. The increasing transport
$T$ from $g_*$ to $g$ is a smooth diffeomorphism. Both $T'$ and $1/T'$
are bounded: on the two tails $g=\epsilon g_*$, and the tail masses
of $g_*$ are asymptotic to constants times $e^{-\lambda |u|}$.
The distribution function identities give $T(u)-u=O(1)$ on each tail
and $T'=g_*/(g\circ T)\to1$ there.

For $z>0$ the function $t\mapsto\log\sinh(e^t)$ is convex, since its
second derivative at $t=\log z$ is
$z(\sinh z\cosh z-z)/\sinh^2 z>0$. For $u<v$, its tangent inequality,
followed by Jensen's inequality, gives
\[
 Q(T(v)-T(u))-Q(v-u)
 \ge (v-u)\coth(v-u)\log\frac{T(v)\!-\!T(u)}{v-u}
 \ge \coth(v-u)\int_u^v\log T'(w)dw .
\]
Define
\[
 A(w)=\lambda\iint_{u<w<v}g_*(u)g_*(v)\coth(v-u)\,du\,dv .
\]
This function is integrable because $(v-u)\coth(v-u)\le1+v-u$ and
$g_*$ has finite first moment. In fact $A=g_*$ almost everywhere.
To see this, transport $g_*$ by $u\mapsto u+\tau h(u)$ for a smooth
compactly supported $h$ and small $\tau$. Its density variation at zero
is $-(g_*h)'$. Stationarity $\log g_*+\lambda Q*g_*=\mathrm{constant}$
therefore gives zero first variation of $\mathcal F$. Differentiating the
same energy in transport coordinates gives
\[
 \int g_* h'
 =\lambda\iint_{u<v}g_*(u)g_*(v)\coth(v-u)(h(v)-h(u))\,du\,dv
 =\int A h'.
\]
The logarithmic singularity of $Q$ is locally integrable and its growth
at infinity is linear, justifying differentiation in density coordinates.
In transport coordinates the factor $h(v)-h(u)$ cancels the pole of
$\coth(v-u)$ at the diagonal. Thus $A-g_*$ is constant as a distribution,
and its integrability makes that constant zero.

The entropy increment under $T$ is $-\int g_*\log T'$. Integrating the
tangent inequality bounds the interaction increment from below by
$\int A\log T'$. The two expressions cancel and give
$\mathcal F(g)\ge\mathcal F(g_*)$.

For an even bounded compactly supported competitor $f$, its transformed
density $g$ is bounded, vanishes for large positive $u$ and is $O(e^u)$
at $-\infty$. It has finite absolute entropy and first moment. Moreover
$\iint |Q(u-v)|g(u)g(v)\,du\,dv<\infty$ by boundedness of $g$, the
integrability of the logarithmic singularity and linear growth of $Q$.
Cut $g$ to $[-m,m]$ and normalize; the energy converges by dominated
convergence. Each bounded compactly supported truncation can be mollified:
the densities stay bounded on a common compact interval, and both entropy
and the logarithmic interaction converge, by dominated convergence along
an almost everywhere convergent subsequence. Mixing each smooth density
with $\epsilon g_*$ and sending $\epsilon$ to zero also preserves the
limit of both terms, since the result is $\epsilon g_*$ outside a
fixed compact set. Thus the inequality holds on the original class.
Conversely, truncating $g_*$ on $[-m,m]$ and normalizing yields bounded
compactly supported competitors in the $x$ variable, with energies
converging to $\mathcal F(g_*)$. This proves sharpness. For any probability
density with all three defining integrals absolutely convergent (omitting
the weight term at $\lambda=1$), one can truncate to
$\{|x|\le m,\ f(x)\le m\}$ and normalize. Dominated convergence in those
integrals extends the same upper bound to that class, which contains $f_*$.

\subsubsection*{Value}
By \eqref{eq:weighted-stationary},
$B_\lambda(f_*)=-\log(\sin a/\pi)+\lambda\int P f_*$. In
$t=|x|^\lambda$ the probability is $(\sin a/a)(1+2t\cos a+t^2)^{-1}dt$.
Substituting $t=\sin(a-\theta)/\sin\theta$ gives the uniform probability
$d\theta/a$ on $(0,a)$, and $1+2t\cos a+t^2=\sin^2 a/\sin^2\theta$.
This proves Theorem \ref{t:loghls}, and Section \ref{sec:hlevel} yields
\[
 C_\lambda^0=1+\lambda\,2^{\lambda-1}\pi
       \exp\!\left(-\frac1a\int_0^a\log\sin\theta\,d\theta\right).
\]
This completes the proof of Theorem \ref{t:level}.

\subsection{Critical attainment on intervals}\label{sec:hattain}
In this subsection we prove Theorem \ref{t:attain}.

\subsubsection*{The centered interval}
Fix $\lambda$ and put $a=\pi\lambda/2$. Write $F=f_*$ for the optimizing
probability density
\[
 F(x)=\frac{\sin a}{\pi}\,
       \frac{|x|^{\lambda-1}}
        {1+2\cos(a)|x|^\lambda+|x|^{2\lambda}},
 \qquad B_*=B_\lambda(F).
\]
Choose $R=\varepsilon^{-1/4}$ and let $f_R$ be the normalized
restriction of $F$ to $R^{-1}<|x|<R$. Then
\begin{equation}
 B_\lambda(f_R)=B_*+O_\lambda(R^{-\lambda}\log R).
 \label{eq:P16}
\end{equation}
To check this, the discarded mass is $O(R^{-\lambda})$ because
$F(x)\asymp|x|^{\lambda-1}$ near zero and
$F(x)\asymp |x|^{-1-\lambda}$ near infinity. The discarded absolute
entropy and logarithmic moment are $O(R^{-\lambda}\log R)$.
The same bound holds for the change of the logarithmic interaction.
Indeed
\[
 \int |\log|x-y|| F(y)\,dy\le C_\lambda(1+\log(1+|x|)).
\]
For the positive part use $|x-y|\le(1+|x|)(1+|y|)$. The negative
part is uniformly bounded by H\"older's inequality, since
$F\in L^{1+\eta}(\mathbb R)$ for some $\eta>0$ and
$\mathbf 1_{|z|<1}|\log|z||$ has finite moments of every finite order.
Renormalization adds an $O(R^{-\lambda})$ error to these quantities.
This proves \eqref{eq:P16}.

As in the proof of the lower bound in Section \ref{sec:hlevel}, let
$f_\varepsilon(x)=\varepsilon^{-1}f_R(x/\varepsilon)$,
$A_\varepsilon=\operatorname{supp} f_\varepsilon$,
$m_\varepsilon=\langle f_\varepsilon,Gf_\varepsilon\rangle$,
and $v_\varepsilon=Gf_\varepsilon/\sqrt{m_\varepsilon}$.
The numerator of the Green kernel on $A_\varepsilon^2$ equals
$2+O((\varepsilon R)^2)$ uniformly. The interaction bound used in
\eqref{eq:P16} therefore gives $\pi m_\varepsilon=\log(1/\varepsilon)+O_\lambda(1)$.
Repeating the Jensen lower bound of Section \ref{sec:hlevel} with this uniform
remainder gives, for $d\mu_\lambda=(\lambda/2)|x|^{\lambda-1}dx$,
\[
 \int_{A_\varepsilon} e^{\pi\lambda v_\varepsilon^2}\,d\mu_\lambda
 \ge C_\lambda^0-1-
    O_\lambda\big(R^{-\lambda}\log R+(\varepsilon R)^2\big).
\]
On the complement, $\exp(\pi\lambda v_\varepsilon^2)\ge
1+\pi\lambda v_\varepsilon^2$. For $x\in[1/3,1/2]$ and
$y\in A_\varepsilon$, the formula for $G(x,y)$ gives
$G(x,y)\ge c>0$ for all sufficiently small $\varepsilon$.
Hence $v_\varepsilon^2\ge c^2/m_\varepsilon$ on this fixed interval,
which lies outside $A_\varepsilon$. Using
$\mu_\lambda(A_\varepsilon)\le(\varepsilon R)^\lambda$, we obtain
\[
 \int_I e^{\pi\lambda v_\varepsilon^2}\,d\mu_\lambda
 \ge C_\lambda^0+\frac{c_\lambda}{m_\varepsilon}
     -O_\lambda\big((\varepsilon R)^\lambda+
               (\varepsilon R)^2+R^{-\lambda}\log R\big)
 > C_\lambda^0
\]
for $\varepsilon$ small. The final inequality holds because all errors
are $o(1/\log(1/\varepsilon))$. These trials have unit Hilbert energy.
Thus the supremum is strictly above the concentration level. A
rearranged maximizing sequence has a nonzero weak limit, and the
energy splitting and uniform integrability argument of Section \ref{sec:attainment} gives a maximizer. In this Hilbert case, the energy of
the weak remainder has upper limit at most one minus the energy of
the limit. Approximating the limit in energy by a bounded smooth
function leaves a strict gap in that bound, so the weighted endpoint
estimate gives a uniform $L^{1+\delta}(d\mu_\lambda)$ bound for the
exponentials for some $\delta>0$. Vitali's theorem then gives
convergence of the maximizing integrals.

\subsubsection*{Intervals with an off-center pole}
Let $I'=(l,r)$ with $l<0<r$, and $0<\lambda<1$. On this interval
define the concentration excess as the supremum of the limsups of
$\int_{I'}(e^{\pi\lambda v^2}-1)|x|^{\lambda-1}dx$ over
\emph{all} admissible weakly null sequences. Its exact value is
\[
 L_{I'}=c^\lambda\frac{2}{\lambda}(C_\lambda^0-1),
 \qquad c=\frac{-2lr}{r-l}.
\]
The concentration level for the exponential itself with the weight
normalized to be a probability is
$1+L_{I'}/\int_{I'}|x|^{\lambda-1}dx$. The critical supremum is attained.

Indeed, put $d=-(l+r)/(r-l)$ and $x=c t/(1+dt)$. This real fractional
linear map sends $(-1,1)$ to $I'$ and fixes the pole. Its Jacobian
satisfies $dx\,dy/|x-y|^2=dt\,ds/|t-s|^2$ when
$x=c t/(1+dt)$ and $y=c s/(1+ds)$. Thus pullback is an energy
isometry of the zero-exterior spaces, and preserves weakly null
sequences. Their ordinary $L^2$ norms on the two bounded intervals are
equivalent because the Jacobian and its reciprocal are bounded there.
The transformed weight is
$c^\lambda |t|^{\lambda-1}w(t)dt$ with
$w(t)=(1+dt)^{-1-\lambda}$ smooth, positive, and $w(0)=1$.

Let first $0<\lambda<1$. On the centered interval, symmetric decreasing rearrangement bounds
the concentration excess of arbitrary weakly null sequences by
$(2/\lambda)(C_\lambda^0-1)$. Indeed rearrangement decreases the
energy, preserves the vanishing $L^2$ norm, and increases the
weighted exponential excess. Radial sequences give the reverse
inequality for the supremum. Along any admissible weakly null
sequence one has $e^{\pi\lambda v_j^2}-1\to0$ in ordinary $L^1(I)$:
compactness gives convergence in measure, and the unweighted
endpoint gives uniform integrability by choosing $q>1$ with
$q\lambda\le1$. Outside $|t|<\delta$ this also proves convergence
for either weight. On $|t|<\delta$ their ratio after factoring out
$c^\lambda$ differs from one by $O(\delta)$. The centered weighted
endpoint bounds the total exponential excess uniformly, so sending
$\delta\downarrow0$ gives the stated value of $L_{I'}$.

For the trial just constructed, $w=1+O(\varepsilon R)$ on the support of
$f_\varepsilon$. The centered weighted endpoint bounds the core integral,
so this costs $O(\varepsilon R)$ after factoring out $c^\lambda$.
On $[1/3,1/2]$ the factor $w$ has a positive lower bound, preserving
the gain of order $1/m_\varepsilon$. Thus the supremum on $I'$ strictly
exceeds its concentration level. The same energy splitting applies to
arbitrary maximizing sequences, since any remaining weighted endpoint
bound on the centered interval stays valid after multiplication by the
bounded positive function $w$. This gives attainment on $I'$. At
$\lambda=1$ the weight is constant, and the affine map from $(-1,1)$ onto a bounded interval preserves the energy, the zero exterior condition and normalized Lebesgue measure. Hence it preserves the concentration level and the attainment, and the normalized concentration level is $1+2\pi$ on every bounded interval.

\medskip

\noindent\textbf{Acknowledgments}. The second author would like to thank her supervisors, Professors Nguyen Lam and Deping Ye, for their guidance and support.

\medskip

\noindent\textbf{AI Declaration}.
During the preparation of this work the authors used AI tools to draft and revise parts of
the text and proofs, to search and check the literature, and to look for errors in the
arguments. The authors checked every proof and every citation, edited the content as needed,
and take full responsibility for the content of the publication.

\medskip

\noindent\textbf{Conflict of interest}. The authors state no conflict of interest.

\medskip

\noindent\textbf{Data availability}. Not applicable.

\end{document}